\documentclass[a4paper]{article}
\usepackage[latin1]{inputenc}
\usepackage[T1]{fontenc}
\usepackage[francais,english]{babel}
\usepackage{amsmath}
\usepackage{amssymb}

\usepackage{pstricks-add}  
\usepackage{float}     

\usepackage{theorem}
\newtheorem{thm}{Theorem}
\newtheorem{Prop}{Proposition}
\newtheorem{Def}{Definition}
\newtheorem{Lem}{Lemma}

\newtheorem{rk}{Remark}

\newcommand{\R}{\mathbb{R}}
\newcommand{\C}{\mathbb{C}}

\newcommand{\JJ}{\mathcal{J}}

\newcommand{\PP}{\mathcal{P}}
\newcommand{\QQ}{\mathcal{Q}}

\newcommand{\VV}{\mathcal{V}}

\newcommand{\N}{\mathbb{N}}
\newcommand{\lag}{\langle}
\newcommand{\rag}{\rangle}
\newcommand{\Tt}{T_{\theta}}

\newcommand{\vpm}{v_{\pm}}
\newcommand{\wpm}{w_{\pm}}

\numberwithin{equation}{section}

\author{
Karine \textsc{Beauchard}\footnote{CMLS, Ecole Polytechnique, 91 128 Palaiseau cedex, France,
email: Karine.Beauchard@math.polytechnique.fr},
Morgan \textsc{Morancey}\footnote{CMLA, ENS Cachan, CNRS, Universud,
61 avenue du Pr\'{e}sident Wilson, F-94230 Cachan, France,
email: Morgan.Morancey@cmla.ens-cachan.fr}
\thanks{The authors were partially supported by the ``Agence Nationale de la Recherche'' (ANR),
Projet Blanc EMAQS number ANR-2011-BS01-017-01}
}

\title{Local controllability of 1D Schrödinger equations with bilinear control and minimal time}
\date{}

\begin{document}
\maketitle

\begin{abstract}
We consider a linear Schrödinger equation, on a bounded interval, with bilinear control.

In \cite{KB-CL}, Beauchard and Laurent prove that, under an appropriate non degeneracy assumption, 
this system is controllable, locally around the ground state, in arbitrary time. 
In \cite{JMC-CRAS-Tmin}, Coron proves that a positive minimal time is required for this controllability result,
on a particular degenerate example.

In this article, we propose a general context for the local controllability to hold in large time,
but not in small time. The existence of a positive minimal time is closely related to the behaviour
of the second order term, in the power series expansion of the solution.
\end{abstract}

\section{Introduction}

\subsection{The problem}

Let us consider the 1D Schrödinger equation
\begin{equation} \label{Schro:syst}
\left\lbrace \begin{array}{ll}
i \partial_t \psi(t,x) = - \partial_x^2 \psi(t,x) - u(t) \mu(x) \psi(t,x), &  (t,x) \in \mathbb{R}\times (0,1),\\
\psi(t,0)=\psi(t,1)=0,                                                     &   t \in \mathbb{R}.
\end{array}\right.
\end{equation}
Such an equation arises in the modelization of a quantum particle, in an infinite
square potential well, in a uniform electric field with amplitude $u(t)$.
The function $\mu:(0,1) \rightarrow \mathbb{R}$ is the dipolar moment of the particle.
The system (\ref{Schro:syst}) is a bilinear control system in which
the state is the wave function $\psi$, with $\|\psi(t)\|_{L^2(0,1)}=1, \forall t \in \mathbb{R}$
and the control is the real valued function $u$. 
\\

In this article, we study the minimal time required for
the local controllability of (\ref{Schro:syst}) around the ground state.
Before going into details, let us introduce several notations.
The operator $A$ is defined by
\begin{equation} \label{def:A}
\begin{array}{cc}
D(A):=H^2 \cap H^1_0((0,1),\mathbb{C}), 
&
A \varphi := - \frac{d^2 \varphi}{dx^2}.
\end{array}
\end{equation}
Its eigenvalues and eigenvectors are
\begin{equation} \label{vap_vep}
\lambda_k:=(k\pi)^2, \quad 
\varphi_k(x):=\sqrt{2} \sin (k \pi x), \forall k \in \mathbb{N}^*.
\end{equation}
The family $(\varphi_k)_{k \in \mathbb{N}^*}$ is an orthonormal basis
of $L^2((0,1),\mathbb{C})$ and
$$\psi_k(t,x) := \varphi_k(x) e^{-i\lambda_k t}, \forall k \in \mathbb{N}^*$$
is a solution of (\ref{Schro:syst}) with $u \equiv 0$ called eigenstate,
or ground state, when $k =1$. We denote by $\mathcal{S}$ the unit $L^2((0,1),\mathbb{C})$-sphere.
In this article, we consider two types of initial conditions for (\ref{Schro:syst}):
the ground state
\begin{equation} \label{IC}
\psi(0,x)=\varphi_1(x), \quad x \in (0,1),
\end{equation}
or an arbitrary one
\begin{equation} \label{IC_qlq}
\psi(0,x)=\psi_0(x), \quad x \in (0,1).
\end{equation}
\\

Now, let us define the concept of local controllability used in this article.
\begin{Def}
Let $T>0$, $X$ and $Y$ be normed spaces such that $X \subset L^2((0,1),\mathbb{C})$ and $Y \subset L^2((0,T),\mathbb{R})$. 
The system (\ref{Schro:syst}) is 
\textbf{controllable in $X$, locally around the ground state, with controls in $Y$, in time $T$}, 
if, for every $\epsilon>0$, there exists $\delta>0$ such that,
for every $\psi_f \in \mathcal{S} \cap X$ with $\| \psi_f - \psi_1(T) \|_{X} < \delta$,
there exists $u \in Y$ with $\|u\|_Y < \epsilon$ such that
the solution of the Cauchy problem (\ref{Schro:syst})-(\ref{IC}) satisfies $\psi(T)=\psi_f$.
\end{Def}
In particular, this definition requires that arbitrarily small motions may be done with arbitrarily small controls.

\noindent
In this introduction, we first recall two previous results concerning local controllability of systems similar to  (\ref{Schro:syst}). We present a positive result in arbitrary time and a setting for which there exists a positive minimal time. Then, we present the main results of this article i.e. we give a precise setting where local controllability hold in time larger than a minimal time and fails otherwise. We end by a short bibliography and by setting some notations.

\subsection{A first previous result}

First, let us introduce the normed spaces
\begin{equation} \label{def:Hs}
H^s_{(0)}((0,1),\mathbb{C}):=D(A^{s/2}), 
\quad 
\| \varphi \|_{H^s_{(0)}} := \left( \sum\limits_{k=1}^{\infty}
| k^s \langle \varphi , \varphi_k \rangle |^2 \right)^{1/2},
\quad \forall s>0.
\end{equation}
The following result, proved in \cite{KB-CL}, emphasizes that the local controllability holds
in any positive time when the dipolar moment $\mu$ satisfies an appropriate non-degeneracy assumption.

\begin{thm} \label{Thm:KB-CL}
Let $T>0$ and $\mu \in H^{3}((0,1),\mathbb{R})$ be such that 
\begin{equation} \label{hyp_mu}
\exists c>0 \text{ such that } 
\frac{c}{k^3} \leqslant |\langle \mu \varphi_1 , \varphi_k \rangle |, 
\forall k \in \mathbb{N}^*.
\end{equation}
There exists $\delta>0$ and a $C^1$ map 
$\Gamma: \Omega_T  \rightarrow  L^2((0,T),\mathbb{R})$
where
$$\Omega_T := \{ \psi_f \in \mathcal{S} \cap H^3_{(0)}((0,1),\mathbb{C}) ;
\| \psi_f - \psi_1(T) \|_{H^3} < \delta  \},$$
such that, $\Gamma( \psi_1(T) ) =0$ and for every $\psi_f \in  \Omega_T$,
the solution of the Cauchy problem (\ref{Schro:syst})-(\ref{IC})
with control $u:=\Gamma(\psi_f)$ satisfies $\psi(T)=\psi_f$.
\end{thm}

First, let us remark that the assumption (\ref{hyp_mu}) holds for example with $\mu(x)=x^2$.
Actually, it holds generically in $H^3((0,1),\mathbb{R})$ (see \cite[Proposition 16]{KB-CL}).
Indeed, for $\mu \in H^3((0,1),\mathbb{R})$,  three integrations by part and the Riemann-Lebesgue Lemma prove that
\begin{equation} \label{3IPP+RL}
\langle \mu \varphi_1,\varphi_k \rangle
= 2 \int\limits_0^1 \mu(x)\sin(\pi x)\sin( k \pi x) dx =
\frac{4[(-1)^{k+1}\mu'(1)-\mu'(0)]}{k^3 \pi^2} + \underset{k \rightarrow + \infty}{o} \left( \frac{1}{k^3} \right).
\end{equation}
In particular, a necessary (but not sufficient) condition on $\mu$ for (\ref{hyp_mu}) to be satisfied is
$\mu'(1)\pm\mu'(0)\neq 0$.
\\

Note that the function spaces in Theorem \ref{Thm:KB-CL} are optimal. Indeed, they are the same as for the well posedness
of the Cauchy problem (\ref{Schro:syst})(\ref{IC}) (see Proposition \ref{WP-CYpb}).
\\

Finally, let us summarize the proof of Theorem \ref{Thm:KB-CL} in \cite{KB-CL}.
This proof relies on the linear test (see \cite[Chapter 3.1]{JMC-book}), 
the inverse mapping theorem and a regularizing effect.
In particular, the assumption (\ref{hyp_mu}) is necessary for the linearized system to be controllable in
$H^3_{(0)}((0,1),\mathbb{C})$ with controls in $L^2((0,T),\mathbb{R})$. When one of the coefficients
$\langle \mu \varphi_1,\varphi_k\rangle$ vanishes, then the linearized system is not controllable anymore
and the strategy of \cite{KB-CL} fails.

\subsection{A second previous result}

The first article in which a positive minimal time is proved, 
for the local controllability of systems similar to (\ref{Schro:syst}), is \cite{JMC-CRAS-Tmin}.
In this reference, Coron considers the control system
\begin{equation} \label{Schro:syst+S+D}
\left\lbrace \begin{array}{ll}
i \partial_t \psi(t,x) = - \partial_x^2 \psi(t,x) - u(t) (x-1/2) \psi(t,x), &  (t,x) \in \mathbb{R}\times(0,1),\\
\psi(t,0)=\psi(t,1)=0,                                                      & t \in \mathbb{R},                \\
s'(t)=u(t), \quad  d'(t)=s(t),                                              & t \in \mathbb{R},                \\
\end{array}\right.
\end{equation}
where the state is $(\psi,s,d)$ and the control is the real valued function $u$.
This system represents a quantum particle in a moving box: 
$u,s,d$  are the acceleration, the speed and the position of the box.

Note that, here, the relation (\ref{hyp_mu}) is not satisfied:
$$\langle (x-1/2) \varphi_1 , \varphi_k \rangle
= \left\lbrace \begin{array}{l}
0 \text{ if } k \text{ is odd,}\\
\frac{8 k}{\pi^2(k^2-1)^2} \text{ if } k \text{ is even,}
\end{array}\right.$$
thus Theorem \ref{Thm:KB-CL} does not apply.

On one hand, it is proved in \cite{KB-JMC} that this system is 
controllable in $H^7_{(0)}((0,1),\mathbb{C}) \times \mathbb{R} \times \mathbb{R}$, 
locally around the ground state $(\psi=\psi_1,s=0,d=0)$,
with controls $u \in L^\infty((0,T),\mathbb{R})$, 
in time $T$ large enough.

On the other hand, Coron proved in \cite{JMC-CRAS-Tmin} that 
this local controllability does not hold in
arbitrary time: contrary to Theorem \ref{Thm:KB-CL}, 
a positive minimal time is required for the local controllability.
Precisely, Coron proved the following statement.

\begin{thm} \label{ThmCoronCRAS}
There exists $\epsilon>0$ such that, 
for every $\overline{d} \neq 0$ and $u \in L^2((0,\epsilon),\mathbb{R})$ satisfying $|u(t)|<\epsilon, \forall t \in (0,\epsilon)$,
the solution $(\psi,s,d) \in C^0([0,\epsilon],H^1_0((0,1),\mathbb{C})) \times C^0([0,\epsilon],\mathbb{R}) \times C^1([0,\epsilon],\mathbb{R})$
of (\ref{Schro:syst+S+D}) such that $(\psi,s,d)(0)=(\psi_1(0),0,0)$ satisfies
$(\psi,s,d)(\epsilon)\neq(\psi_1(\epsilon),0,\overline{d})$.
\end{thm}

The goal of this article is to go further in this analysis: 
\begin{enumerate}
\item we propose a general context for the minimal time to be positive
(in particular, the variables $s$ and $d$ are not required anymore in the state),
\item we propose a sufficient condition for the local controllability 
to hold in large time; this assumption is compatible with the previous context and weaker than (\ref{hyp_mu}),
\item we work in an optimal functional frame, for instance, 
our non controllability result requires $u$ small in $L^2$-norm, not in $L^\infty$-norm as in Theorem \ref{ThmCoronCRAS},
\item we perform a first step toward the characterization of the minimal time.
\end{enumerate}

\subsection{Main results of this article}

The first result of this article is the following one.

\begin{thm} \label{Main_thm_0}
Let $K \in \mathbb{N}^*$, $\mu \in H^3((0,1),\mathbb{R})$ be such that
\begin{equation} \label{hyp_mu_0}
\langle \mu \varphi_1 , \varphi_K \rangle = 0 \quad \text{ and } \quad A_K:=\langle (\mu')^2 \varphi_1 , \varphi_K \rangle \neq 0,
\end{equation} 
and $\alpha_K \in \{-1,+1\}$ be defined by
\begin{equation} \label{def:alpha}
\alpha_K:=\text{sign} (A_K).
\end{equation}
There exists $T_K^*>0$ such that, for every $T<T_K^*$,
there exists $\epsilon>0$ such that, for every $u \in L^2((0,T),\mathbb{R})$ with 
\begin{equation} \label{upt}
|| u ||_{L^2(0,T)} < \epsilon
\end{equation}
the solution of (\ref{Schro:syst})(\ref{IC}) satisfies 
$\psi(T) \neq [\sqrt{1-\delta^2} \varphi_1 + i \alpha_K \delta \varphi_K ] e^{-i\lambda_1 T}$ for every $\delta>0$.
\end{thm}

First, we remark that the assumption (\ref{hyp_mu_0}) holds, for example, with $\mu(x)=(x-1/2)$ and $K=1$. 
In particular, Theorem \ref{Main_thm_0} applies to the particular case studied by Coron in \cite{JMC-CRAS-Tmin}. 
Thus, the variables $(s,d)$ are not required in the state for the minimal time to be positive. 
Moreover, the control $u$ does not need to be small in $L^\infty(0,T)$ as in Theorem \ref{ThmCoronCRAS}: 
a smallness assumption in $L^2(0,T)$ is sufficient.
\\

Note that the validity of the same result without the assumption '$A_K \neq 0$' is an open problem
(see remark \ref{RK:AK<>0} for technical reasons).
A possible (but not optimal) value of $T_K^*$ is given in (\ref{def:TK*}).
The proof of Theorem \ref{Main_thm_0} relies on an expansion of the solution to the second order.
\\

The second result of this article is the following one.

\begin{thm} \label{thm:control_Tlarge}
Let $\mu \in H^3((0,1),\mathbb{R})$ be such that 
\begin{equation} \label{mu'pm}
\mu'(0)\pm\mu'(1)\neq 0.
\end{equation}
Then, the system (\ref{Schro:syst}) is controllable in $H^3_{(0)}((0,1),\mathbb{C})$,
locally around the ground state, with controls $u \in L^2((0,T),\mathbb{R})$, in large enough time $T$.
\end{thm}

A direct consequence of Theorems \ref{Main_thm_0} and \ref{thm:control_Tlarge} is the following result.

\begin{thm}
Let $\mu \in H^3((0,1),\mathbb{R})$ be such that (\ref{hyp_mu_0}) and (\ref{mu'pm}) hold for some $K \in \mathbb{N}^*$.
Then, there exists $T_{min}>0$ such that the controllability of (\ref{Schro:syst}) in $H^3_{(0)}((0,1),\mathbb{C})$,
locally around the ground state, with controls in $L^2((0,T),\mathbb{R})$ does not hold when $T<T_{min}$, 
and holds when $T>T_{min}$.
\end{thm}

First, we remark that the assumption (\ref{mu'pm}) is weaker than (\ref{hyp_mu})
and that the assumptions  (\ref{hyp_mu_0}) and (\ref{mu'pm}) are compatible: consider, for instance 
$\mu(x):=x^2-\langle x^2\varphi_1,\varphi_2\rangle \varphi_2/\varphi_1$.
\\

Note that an explicit upper bound $T_\sharp$ for the minimal time $T_{min}$ is proposed in 
the proof (see (\ref{def_Tmin})).
\\

We emphasize that, when $\mu'(0)=\mu'(1)=0$, then, the appropriate functional frame stops to be
 $(\psi \in H^3_{(0)},u \in L^2)$. For instance, with the tools developed in this article,
one may prove: if $L \in \mathbb{N}$, $\mu \in H^{2L+3}((0,1),\mathbb{R})$ are such that 
$\mu^{(2k+1)}(0)=\mu^{(2k+1)}(1)=0$ for $k=0,...,L-1$ and 
$\mu^{(2L+1)}(0) \pm \mu^{(2L+1)}(1) \neq 0$,
then, the system (\ref{Schro:syst}) is controllable in $H^{2L+3}_{(0)}((0,1),\mathbb{C})$,
locally around the ground state, 
with controls in $L^2((0,T),\mathbb{R})$, 
in large enough time $T$.
\\

Finally, we summarize the proof of Theorem \ref{thm:control_Tlarge}.
Under assumption (\ref{mu'pm}), only a finite number of the coefficients $\langle\mu\varphi_1,\varphi_k\rangle$
vanish (see (\ref{3IPP+RL})). Thus, the linearized system around the ground state is not controllable 
along  a finite number of directions. We will see that all of these directions are recovered at the second order.
Moreover, all these directions excepted one, present a rotation phenomena in the complex plane,
for the null input solution. This idea of using a power series expansion and exploiting a rotation phenomena 
was first used on a Korteweg-de Vries equation by Cerpa and Crépeau in \cite{Crepeau-Cerpa}.
However, their strategy has to be adapted in our situation, 
because one lost direction does not exhibit a rotation phenomenon (see Remark \ref{rk_rotation}).
\\

Under a weaker assumption than (\ref{mu'pm}) and still in the framework $(\psi \in H^3_{(0)},u \in L^2)$,
we prove the following result.

\begin{thm} \label{thm:control_partiel_Tlarge_N}
Let $\mu \in H^3((0,1),\mathbb{R})$ be such that
\begin{equation} \label{mu'_partiel}
\mu'(0) = \mu'(1)\neq 0 \quad (\text{resp.} \mu'(0) = -\mu'(1)\neq 0).
\end{equation}
For $N \in \mathbb{N}^*$, we define $\mathcal{N}_N := \left\{ k \in \N^* \, ; \, k \text{ is odd and } k \leq N \text{ or } k \text{ is even } \right\}$ (resp. $\mathcal{N}_N := \left\{ k \in \N^* \, ; \, k \text{ is even and } k \leq N \text{ or } k \text{ is odd } \right\}$).
Let $\mathbb{P}_N$ be the orthogonal projection from $L^2((0,1),\mathbb{C})$ to 
$\mathbb{V}_N:=\text{Span}\{ \varphi_k ; k \in \mathcal{N}_N \}$.
Then, for every $\epsilon>0$, there exists $T>0$ and $\delta>0$ such that,
for every $\widetilde{\psi}_f \in \mathbb{V}_N \cap H^3_{(0)}(0,1)$ with
$\| \widetilde{\psi}_f - \mathbb{P}_N \psi_1(T) \|<\delta$,
there exists $u \in L^2(0,T)$ with $\|u\|_{L^2}<\epsilon$ such that
the solution of (\ref{Schro:syst})-(\ref{IC}) satisfies
$\mathbb{P}_N \psi(T)=\widetilde{\psi}_f$.
\end{thm}

The sketch of the proof is the following. Under assumption (\ref{mu'_partiel}), we prove that
\begin{enumerate}
\item an infinite number of directions are controlled at the first order, in any positive time,
\item all the lost directions are recovered either at the second order, or at the third order,
\item any direction corresponding to vanishing first and second orders,
are recovered at the third order in arbitrary time.
\end{enumerate}
Note that even if $\mu'(0)=\mu'(1) \neq 0$ (resp.  $\mu'(0)=-\mu'(1) \neq 0$),
one may sometimes control the whole wave function $\psi$ in large time.
For instance in \cite{KB-JMPA}, the local controllability in $H^7_{(0)}((0,1),\mathbb{C})$, 
with controls in $H^1_0((0,T),\mathbb{R})$, in large time $T$, 
is proved for $\mu(x)=(x-1/2)$, with the return method.

\subsection{A review about control of bilinear systems}

The first controllability result for bilinear Schrödinger equations such as (\ref{Schro:syst}) 
is negative and proved by Turinici \cite{Turinici},
as a corollary of a more general result by Ball, Marsden and Slemrod 
\cite{ball-marsden-slemrod}. Then, it has been adapted to nonlinear Schrödinger equations 
in \cite{Teismann-et-al} by Ilner, Lange and Teismann.
Because of such noncontrollability results,
these equations have been considered as non controllable for a long time.
However, progress have been made and this question is now better understood. 
\\

Concerning exact controllability issues, 
local results for 1D models have been proved in \cite{KB-JMPA,SchroLgVar} by Beauchard;
almost global results have been proved in \cite{KB-JMC}, by Coron and Beauchard.
In \cite{KB-CL}, Beauchard and Laurent proposed an important simplification of the above proofs.
In \cite{JMC-CRAS-Tmin}, Coron proved that a positive minimal time may be required for the local controllability of the 1D model.
In \cite{ondes}, Beauchard studied the minimal time for the local controllability of 1D wave equations with bilinear controls. 
In this reference, the origin of the minimal time is the linearized system,
whereas in the present article, the minimal time is related to the nonlinearity of the system.
Exact controllability has also been studied in infinite time by Nersesyan and Nersisian in \cite{VNHN,VNHN2}. 
\\

Now, we quote some approximate controllability results.
In \cite{KB-MM} Mirrahimi and Beauchard proved the global approximate controllability,
in infinite time, for a 1D model and in \cite{MM} Mirrahimi proved a similar result
for equations involving a continuous spectrum.
Approximate controllability, in finite time, has been proved for 
particular models by Boscain and Adami in \cite{Boscain-Adami},
by using adiabatic theory and intersection of the eigenvalues in the space of controls.
Approximate controllability, in finite time, for more general models, have been studied
by three teams, with different tools:
by Boscain, Chambrion, Mason, Sigalotti \cite{Chambrion-et-al, Chambrion-et-al2, Chambrion-et-al3},  
with geometric control methods;
by Nersesyan \cite{Nersesyan1,Nersesyan2}
with feedback controls and variational methods;
and by Ervedoza and Puel \cite{Ervedoza-Puel}
thanks to a simplified model.
\\

Optimal control techniques have also been investigated for Schr\"odinger 
equations with a non linearity of Hartee type
in \cite{Baudouin,Baudouin-Kavian-Puel} by Baudouin, Kavian, Puel
and in \cite{CLB-Cances-Pilot} by Cances, Le Bris, Pilot.
An algorithm for the computation of such optimal controls is studied
in \cite{Baudouin-Salomon} by Baudouin and Salomon.
\\

Finally, we quote some references concerning bilinear wave equations. 
In \cite{KhapalovCOCV,KhapalovDCDS,KhapalovDCDIS},
Khapalov considers nonlinear wave equations with bilinear controls. He proves the global
approximate controllability to nonnegative equilibrium states.

\subsection{Notations}

We introduce some conventions and notations valid in all this article.
Unless otherwise specified, the functions considered are complex valued and,
for example, we write $H^1_0(0,1)$ for $H^1_0((0,1),\mathbb{C})$. When the 
functions considered are real valued, we specify it and we write, for example,
$L^2((0,T),\mathbb{R})$. 
The same letter $C$ denotes a positive constant, that can change from one line
to another one. If $(X,\|.\|)$ is a normed vector space, $x \in X$ and $R>0$, $B_X(x,R)$ denotes
the open ball $\{ y \in X ; \|x-y\| < R \}$ and
$\overline{B}_X(x,R)$ denotes the closed ball $\{ y \in X ; \|x-y\| \leqslant R \}$.
We denote by $\langle \cdot , \cdot \rangle$ the $L^2(0,1)$ hermitian inner product
$$\langle f , g \rangle = \int_0^1 f(x) \overline{g(x)} dx,$$
and by $T_{\mathcal{S}} \varphi := \{ \xi \in L^2(0,1);\Re \langle \varphi,\xi\rangle=0 \}$
the tangent space to $\mathcal{S}$ at any point $\varphi \in \mathcal{S}$.
We also introduce for any $s>0$, the normed spaces
$$ h^s(\N^*,\C) := \left\{ a=(a_k)_{k \in \N^*} \in \C^{\N^*} \, ; \, \sum_{k=1}^{+\infty} |k^s a_k|^2 < + \infty \right\}, \quad ||a||_{h^s} := \left( \sum_{k=1}^{+\infty} |k^s a_k|^2 \right)^{1/2}.$$

\subsection{Structure of this article}

In Section \ref{sec:WP}, we recall a well posedness result concerning system (\ref{Schro:syst}).
In Section \ref{sec:thm_0}, we prove Theorem \ref{Main_thm_0}.
In Section \ref{sec:thm:control_Tlarge}, we prove Theorem \ref{thm:control_Tlarge} 
thanks to power series expansions to the second order as in \cite{Crepeau-Cerpa} (see also (\cite[Chapter 8]{JMC-book})).
In Section \ref{sec:N}, we prove Theorem \ref{thm:control_partiel_Tlarge_N} 
thanks to power series expansions to the order 2 and 3.
In Section \ref{sec:caract}, we perform a first step toward the characterization of the minimal time, 
in a favorable situation. Finally, in Section \ref{sec:ccl}, we gather several concluding remarks and perspectives.

\section{Well posedness}
\label{sec:WP}

This section is dedicated to the well posedness of the Cauchy problem
\begin{equation} \label{CY}
\left\lbrace \begin{array}{ll}
i \partial_t \psi = - \partial_x^{2} \psi - u(t) \mu(x) \psi - f(t,x), &   (t,x) \in(0,T)\times(0,1),\\
\psi(t,0)=\psi(t,1)=0,                                                 & t \in (0,T),\\
\psi(0,x)=\psi_{0}(x),                                                 & x \in (0,1).
\end{array}\right.
\end{equation}
proved in \cite[Proposition 3]{KB-CL}.

\begin{Prop} \label{WP-CYpb}
Let $\mu \in H^{3}((0,1),\mathbb{R})$, 
$T>0$, 
$\psi_0 \in H^3_{(0)}(0,1)$, 
$f\in L^2((0,T),H^3 \cap H^1_0)$
and $u \in L^2((0,T),\mathbb{R})$.
There exists a unique weak solution of (\ref{CY}), i.e.
a function $\psi \in C^0([0,T],H^3_{(0)})$ such that the following equality 
holds in $H^{3}_{(0)}(0,1)$ for every $t \in [0,T]$,
\begin{equation} \label{SolutionFaible}
\psi(t)=e^{-iAt} \psi_{0} +i \int_{0}^{t} e^{-iA(t-\tau)}
[ u(\tau) \mu \psi(\tau) + f(\tau) ] d\tau.
\end{equation}
Moreover, for every $R>0$, there exists $C=C(T,\mu,R)>0$ such that,
if $\|u\|_{L^2(0,T)} < R$, then this weak solution satisfies
\begin{equation} \label{majo}
\| \psi \|_{C^{0}([0,T],H^{3}_{(0)})} \leqslant 
C \Big( \| \psi_{0} \|_{H^{3}_{(0)}} + \|f\|_{L^{2}((0,T),H^{3} \cap H^1_0(0,1))} \Big).
\end{equation}
If $f \equiv 0$ then 
\begin{equation} \label{NL2=1}
\|\psi(t)\|_{L^2(0,1)} = \|\psi_0\|_{L^2(0,1)}, \forall t \in [0,T].
\end{equation}
\end{Prop}

\section{Examples of impossible motions in small time}
\label{sec:thm_0}

The goal of this section is to prove Theorem \ref{Main_thm_0}.

\subsection{Heuristic}
\label{subsec:Heuristic}

Since we are interested in small motions around the trajectory $(\psi=\psi_1,u=0)$, with small controls,
it is natural to try to do them, in a first step, with the first and the second order terms.
We consider a control $u$ of the form $u=0+\epsilon v + \epsilon^2 w$.
Then, formally, the solution $\psi$ of (\ref{Schro:syst})(\ref{IC}) writes
$\psi=\psi_1+\epsilon \Psi + \epsilon^2 \xi + o(\epsilon^2)$ where
\begin{equation} \label{Schro_lin_eq}
\left\lbrace \begin{array}{ll}
i \partial_t \Psi = - \partial_x^2 \Psi - v(t) \mu(x) \psi_1, & (t,x) \in (0,T)\times(0,1), \\
\Psi(t,0)=\Psi(t,1)=0,                                        & t \in (0,T),                \\
\Psi(0,x)=0,                                                  & x \in (0,1),
\end{array}\right.
\end{equation}
\begin{equation} \label{Schro_O2_eq}
\left\lbrace \begin{array}{ll}
i \partial_t \xi = - \partial_x^2 \xi - v(t) \mu(x) \Psi - w(t) \mu(x) \psi_1, &  (t,x) \in (0,T)\times(0,1), \\
\xi(t,0)=\xi(t,1)=0,                                                             & t \in (0,T),                 \\
\xi(0,x)=0,                                                                      & x \in (0,1).          
\end{array}\right.
\end{equation}
From the property $\|\psi(t)\|_{L^2} \equiv 1$, we deduce that
$\Re\langle \Psi(t) , \psi_1(t) \rangle =0$ (i.e. $\Psi(t) \in T_{\mathcal{S}} \psi_1(t)$, $\forall t$) and
\begin{equation}
\label{espace_tangent_ordre2}
\|\Psi(t)\|_{L^2}^2 + 2 \Re \langle \xi(t),\psi_1(t) \rangle \equiv 0.
\end{equation}
We have
\begin{equation} \label{Psi_explicit}
\Psi(T,x)= i \sum\limits_{j=1}^{\infty} 
\langle \mu \varphi_1 , \varphi_j \rangle \int_0^T v(t) e^{i \omega_j t} dt \psi_j(T,x)
\end{equation}
where
\begin{equation} \label{def:omega}
\omega_j:=\lambda_j-\lambda_1, \quad \forall j \in \mathbb{N}^*.
\end{equation}

We assume that (\ref{hyp_mu_0}) holds for some $K \in \mathbb{N}^*$.
By adapting the choice of $v \in L^2((0,T),\mathbb{R})$, $\Psi(T)$ can reach any target in the closed subspace
$\text{Adh}_{H^3_{(0)}(0,1)} [ \text{Span} \{ \varphi_k \, ; \, k \in \JJ \} ]$
where
\begin{equation} \label{def:J}
\JJ:=\{ j \in \mathbb{N}^* ; \langle \mu \varphi_1 , \varphi_j \rangle \neq 0 \}
\end{equation}
(see Proposition \ref{Cor:haraux1} in Appendix); but the complex direction $\langle \Psi(T),\psi_K(T)\rangle$ is lost.
Let us show that, when $T$ is small, the second order term imposes a sign on the component along this lost direction,
preventing the local exact controllability around the ground state.
\\

Using (\ref{Schro_O2_eq}) and (\ref{Psi_explicit}), we get
\begin{equation} \label{CompK_O2}
\langle\xi(T),\psi_K(T)\rangle = Q^2_{K,T}(v),
\end{equation}
where
\begin{equation} \label{def:Q2KT}   
Q^2_{K,T}(v):= \int_0^T v(t) \int_0^t v(\tau) h^2_{K}(t,\tau)  d\tau dt,
\end{equation}
\begin{equation} \label{def:h2K}
h^2_K(t,\tau):= - \sum\limits_{j=1}^\infty \langle\mu\varphi_K,\varphi_j\rangle \langle\mu\varphi_j,\varphi_1\rangle
e^{i[(\lambda_K-\lambda_j)t + (\lambda_j-\lambda_1) \tau ]}.
\end{equation}
The index $2$ in $Q^2_{K,T}$ and $h^2_K$ is related to the fact that $\xi$ is the second order of the power series expansion.
Integrations by part show that
\begin{equation} \label{mu_phi1_mhik}
|\langle \mu \varphi_K , \varphi_j \rangle | \text{ and }
|\langle \mu \varphi_1 , \varphi_j \rangle | \leqslant \frac{C}{j^3}, \quad \forall j \in \mathbb{N}^*,
\end{equation}
for some constant $C=C(\mu)>0$, thus $h^2_K \in C^0(\mathbb{R}^2,\mathbb{C})$ and the quadratic form $Q^2_{K,T}$ is well defined on 
$L^2((0,T),\mathbb{R})$.
In particular,
\begin{equation} \label{comp1_O2}
\Im [ \langle \xi(T),\varphi_K e^{-i\lambda_1 T} \rangle ] 
=  \widetilde{Q}^2_{K,T}(v)
\end{equation}
where
\begin{equation} \label{def:Q2KTtilde}
\widetilde{Q}^2_{K,T}(v):=\int_0^T v(t) \int_0^t v(\tau) \widetilde{h}^2_{K,T}(t,\tau) d\tau dt,
\end{equation}
\begin{equation} \label{def:h}
\widetilde{h}^2_{K,T}(t,\tau):=\sum\limits_{j=1}^{\infty}
\langle \mu \varphi_K , \varphi_j \rangle  \langle \mu \varphi_j , \varphi_1 \rangle 
\sin[(\lambda_j-\lambda_K)t - \omega_j \tau +(\lambda_K - \lambda_1) T].
\end{equation}
\\

Now, we try to move $\epsilon \Psi(T) + \epsilon^2 \xi(T)$ in the direction of $+ i \alpha_K \varphi_K e^{-i\lambda_1 T}$
(see (\ref{def:alpha}) for the definition of $\alpha_K$).
Since $\Psi(T)$ lives in $\text{Adh}_{H^3_{(0)}(0,1)} [ \text{Span} \{ \varphi_k ; k \neq K  \} ]$,
then, necessarily $\Psi(T)=0$, i.e. $v$ belongs to 
\begin{equation} \label{def:VT}
V_T:=\left\{ 
v \in L^2((0,T),\mathbb{R}) ; \int_0^T v(t) e^{i \omega_j t} dt=0, \forall j \in \JJ
\right\}
\end{equation}
and $\xi(T)=i \delta \alpha_K \varphi_K e^{-i\lambda_1 T}$ for some $\delta>0$. 
Thus the sign of $\widetilde{Q}_{K,T}^2(v)$ has to be $\alpha_K$. 
The following two lemmas show that this is not possible when $T$ is small.

\begin{Lem} \label{Lem:correspondanceFQ}
For every $v \in V_T$, we have
$\widetilde{Q}^2_{K,T}(v) = \mathcal{Q}_{K,T}(S)$ 
where $S(t):=\int_0^t v(\tau) d\tau$ and
\begin{equation} \label{def:QT(S)}
\mathcal{Q}_{K,T}(S) := -A_K \int_0^T S(t)^2 \cos[(\lambda_K-\lambda_1)(t-T)]dt
+ \int_0^T S(t) \int_0^t S(\tau) k_{K,T}(t,\tau) d\tau dt,
\end{equation}
\begin{equation} \label{def:kKT}
k_{K,T}(t,\tau):=\sum\limits_{j=1}^\infty 
(\lambda_j-\lambda_K) \omega_j  \langle \mu \varphi_1 , \varphi_j \rangle\langle \mu \varphi_K , \varphi_j \rangle
\sin[(\lambda_j-\lambda_K)t - \omega_j \tau + (\lambda_K-\lambda_1)T ].
\end{equation}
\end{Lem}

\begin{rk}
Note that $\mathcal{Q}_{K,T}$ is well defined on $L^2(0,T)$ because
$k_{K,T} \in L^\infty(\mathbb{R}\times\mathbb{R})$ (see (\ref{mu_phi1_mhik})).
\end{rk}

\noindent \textbf{Proof of Lemma \ref{Lem:correspondanceFQ}:}
Let $T>0$ and $v \in V_{T}-\{0\}$. 
Integrations by parts show that, for every $j \in \JJ$,
\begin{align*}
& \int_0^T v(t) \int_0^t v(\tau) e^{i [ (\lambda_j-\lambda_K)t - \omega_j \tau ]} d\tau dt
\\ = &
- \int_0^T S(t) \left( v(t) e^{i(\lambda_1-\lambda_K)t} 
+ i (\lambda_j-\lambda_K) \int_0^t v(\tau) e^{i [ (\lambda_j-\lambda_K)t - \omega_j \tau ]}  d\tau 
\right) dt
\\ = &
- \frac{1}{2} S(T)^2 e^{i(\lambda_1-\lambda_K)T} + \frac{i(\lambda_1-\lambda_K)}{2} \int_0^T S(t)^2  e^{i(\lambda_1-\lambda_K)t} dt
\\ &
- i (\lambda_j-\lambda_K) \int_0^T S(t) \left(
S(t)e^{i(\lambda_1-\lambda_K)t} +i \omega_j \int_0^t S(\tau) e^{i [ (\lambda_j-\lambda_K)t - \omega_j \tau ]}  d\tau
\right) dt
\\ = &
- \frac{1}{2} S(T)^2 e^{i(\lambda_1-\lambda_K)T}
- i \left( \lambda_j - \frac{\lambda_1+\lambda_K}{2} \right) \int_0^T S(t)^2 e^{i(\lambda_1-\lambda_K)t} dt
\\ &
+(\lambda_j-\lambda_K) \omega_j \int_0^T S(t) \int_0^t S(\tau) e^{i[(\lambda_j-\lambda_K)t-\omega_j \tau]} d\tau dt.
\end{align*}
The relations
$$\sum\limits_{j=1}^{\infty} \langle \mu \varphi_1 , \varphi_j \rangle \langle \mu \varphi_K , \varphi_j \rangle
= \langle \mu \varphi_1 , \mu \varphi_K \rangle,$$
\begin{equation} \label{AK_dvp}
\sum\limits_{j=1}^{\infty} \left( \lambda_j - \frac{\lambda_1+\lambda_K}{2} \right) 
\langle \mu \varphi_1 , \varphi_j \rangle \langle \mu \varphi_K , \varphi_j \rangle
= \langle (\mu')^2 \varphi_1 , \varphi_K \rangle = A_K.
\end{equation}
give the conclusion. \hfill  $\blacksquare$

\begin{Lem} \label{Lem:T_2}
Let $\mu \in H^3((0,1),\mathbb{R})$ be such that (\ref{hyp_mu_0}) holds for some $K \in \mathbb{N}^*$. 
There exists $T_K^*>0$ such that, for every $T<T_K^*$ 
\begin{equation} \label{QKT_coercive}
\mathcal{Q}_{K,T}(S)
\left\lbrace \begin{aligned}
\leqslant - \frac{A_K}{4} \int_0^T S(t)^2 dt  \text{ if } A_K>0,\\
\geqslant - \frac{A_K}{4} \int_0^T S(t)^2 dt  \text{ if } A_K<0,
\end{aligned}\right\},
\forall S \in L^2((0,T),\mathbb{R}).
\end{equation}
\end{Lem}

\begin{rk} \label{RK:AK<>0}
This statement enlightens the importance of the assumption $A_K \neq 0$ in Theorem \ref{Main_thm_0}.
Indeed, if $A_K$ vanishes then we do not know whether the quadratic form $\widetilde{Q}^2_{K,T}$ has a sign on $V_T$
in small time $T$. Note that another integration by parts
(leading to a quadratic form in $\sigma(t):=\int_0^t S$) is not possible,
because of problems of divergence in infinite sums.
\end{rk}

\noindent \textbf{Proof of Lemma \ref{Lem:T_2}:} One may assume that $A_K>0$, $\alpha_K=1$.
We define the quantity
$$C_K:=\sum\limits_{j=1}^\infty 
|(\lambda_j-\lambda_K) \omega_j  \langle \mu \varphi_1 , \varphi_j \rangle\langle \mu \varphi_K , \varphi_j \rangle|.$$
By (\ref{AK_dvp}) and (\ref{hyp_mu_0}), there exists $j \in \mathbb{N}^*-\{1,K\}$ such that
$ \langle \mu \varphi_1 , \varphi_j \rangle\langle \mu \varphi_K , \varphi_j \rangle \neq 0$.
Thus, $C_K>0$. We introduce 
\begin{equation} \label{def:TK*}
T_K^* := \left\lbrace \begin{array}{l}
\frac{|A_1|}{2C_1} \text{ if } K=1, \\
\min\left\{ \frac{|A_K|}{2C_K};\frac{\pi}{3(\lambda_K-\lambda_1)} \right\} \text{ if } K \geqslant 2.
\end{array}\right. 
\end{equation}
Let $T \in (0,T_K^*)$. Using the inequality
$$\cos[(\lambda_K-\lambda_1)(t-T)] \geqslant \frac{1}{2},\quad \forall t \in (0,T),$$
(\ref{def:QT(S)}), (\ref{def:kKT}) and Cauchy-Schwarz inequality we get, for every $S \in L^2((0,T),\mathbb{R})$,
\begin{align*}
\mathcal{Q}_{K,T}(S) 
& \leqslant - \frac{A_K}{2} \int_0^T S(t)^2 dt + C_K \int_0^T |S(t)| \int_0^t |S(\tau)| d\tau dt
\\ & \leqslant
-\frac{1}{2} \left[ A_K - T C_K  \right] \int_0^T S(t)^2 dt.  
\end{align*}
\hfill $\blacksquare$

With additional arguments, one may prove that, for $T<T^*_K$, 
$$\sup\{ \widetilde{Q}^2_{K,T}(v) ; v \in V_T , \|v\|_{L^2}=1 \}=0.$$
The non existence of a  positive constant $c(T)>0$ such that
$$\widetilde{Q}^2_{K,T}(v) \leqslant - c(T) \|v\|_{L^2}^2, \forall v \in V_T, \forall T<T^*_K$$
prevents from proving the non controllability in a simple way.
Our solution relies on the fact that, for $T$ small, the quadratic form $\widetilde{Q}^2_{K,T}$ is
coercive in $S(t):=\int_0^t v(\tau) d \tau$ (see (\ref{QKT_coercive})). 
This justifies several technical developments and the use of an auxiliary system in the next section.

\subsection{Auxiliary system}

We consider the function $\widetilde{\psi}(t,x)$ defined by
\begin{equation} \label{def:psitilde}
\psi(t,x)= \widetilde{\psi}(t,x) e^{i s(t) \mu(x)}
\text{ where } s(t):=\int_0^t u(\tau) d\tau,
\end{equation}
which is a weak solution of
\begin{equation} \label{Schro:syst_aux}
\left\lbrace \begin{array}{ll}
i \partial_t \widetilde{\psi} = - \partial_x^2 \widetilde{\psi}  
-i s(t)[2\mu'(x) \partial_x  \widetilde{\psi} + \mu''(x)  \widetilde{\psi}] + s(t)^2 \mu'(x)^2  \widetilde{\psi}, 
                                                                  & \quad (t,x) \in (0,T)\times(0,1),\\
\widetilde{\psi}(t,0)= \widetilde{\psi}(t,1)=0,                   & \quad t \in (0,T), \\
\widetilde{\psi}(0,x)= \varphi_1(x),                                 & \quad x \in (0,1).
\end{array}\right.
\end{equation}
We deduce from (\ref{def:psitilde}) and Proposition \ref{WP-CYpb} (applied to (\ref{Schro:syst})(\ref{IC})) the following well posedness result for (\ref{Schro:syst_aux}).

\begin{Prop} \label{WP-CYpb_aux}
Let $\mu \in H^{3}((0,1),\mathbb{R})$, $T>0$,  $s \in H^1((0,T),\mathbb{R})$ with $s(0)=0$.
There exists a unique weak solution  $\widetilde{\psi} \in C^0([0,T],H^3 \cap H^1_0(0,1))$ of (\ref{Schro:syst_aux}). 
Moreover, for every $R>0$, there exists $C=C(T,\mu,R)>0$ such that,
if $\|\dot{s}\|_{L^2(0,T)} < R$, then this weak solution satisfies
\begin{equation} \label{Continuite_aux}
\| \widetilde{\psi} \|_{L^{\infty}((0,T),H^3 \cap H^1_0)} \leqslant C.
\end{equation}
\end{Prop}

The proof of Theorem \ref{Main_thm_0} is a direct consequence of the following result.

\begin{thm} \label{Thm:aux}
Let $K \in \mathbb{N}^*$, $\mu \in H^3((0,1),\mathbb{R})$ be such that (\ref{hyp_mu_0}) holds 
and $T_K^*$ be as in Lemma \ref{Lem:T_2}. For every $T<T^*_K$,
there exists $\epsilon>0$ such that
for every $s \in H^1((0,T),\mathbb{R})$ with $s(0)=0$ and $\| \dot{s}\|_{L^2} < \epsilon$,
the solution of the Cauchy problem (\ref{Schro:syst_aux}) satisfies
\begin{equation} \label{CF_aux}
\widetilde{\psi}(T,.) \neq ( \sqrt{1-\delta^2}\varphi_1 + i \alpha_K \delta \varphi_K ) e^{-i\lambda_1 T} e^{i \theta \mu}, \;
\forall \delta>0, \forall \theta \in \mathbb{R}.
\end{equation}
\end{thm}

The proof of Theorem \ref{Thm:aux} requires several steps, 
thus, it is developed in Section \ref{subsec:Proof_thm_aux}.

\subsection{Proof of Theorem \ref{Main_thm_0} thanks to Theorem \ref{Thm:aux}}

Let $T<T^*_K$. Let $\epsilon>0$ be as in Theorem \ref{Thm:aux}.
Let $u \in L^2((0,T),\mathbb{R})$ be such that $||u||_{L^2}< \varepsilon$.
We assume that the solution of the Cauchy problem (\ref{Schro:syst})(\ref{IC}) satisfies
$\psi(T) = (\sqrt{1-\delta^2}\varphi_1 + i \alpha_K \delta \varphi_K) e^{-i\lambda_1 T}$ for some $\delta>0$.
Then, the function $\widetilde{\psi}$ defined by (\ref{def:psitilde}) solves (\ref{Schro:syst_aux}) and satisfies
$\widetilde{\psi}(T) = (\sqrt{1-\delta^2}\varphi_1 + i \alpha_K \delta \varphi_K) e^{-i\lambda_1 T} e^{-is(T)\mu}$.
By Theorem \ref{Thm:aux}, this is impossible.

\subsection{Proof of Theorem \ref{Thm:aux}}
\label{subsec:Proof_thm_aux}

The proof of Theorem \ref{Thm:aux} requires the following preliminary result.

\begin{Prop} \label{Prop:approx_Comp1k}
Let $T>0$, $K \in \mathbb{N}^*$, $\mu \in H^3((0,1),\mathbb{R})$ be such that $\langle \mu \varphi_1,\varphi_K \rangle=0$.
\begin{equation} \label{approx_comp1_O2}
\Big| \Im \langle \widetilde{\psi}(T) , \varphi_K e^{-i\lambda_1 T} \rangle - \QQ_{K,T}(s) \Big| 
= o(||s||_{L^2}^2) \text{ when } ||u||_{L^2} \to 0,
\end{equation}
\begin{equation} \label{approx_O1}
| \Im \langle \widetilde{\psi}(T) , \psi_1(T) \rangle | = o(||s||_{L^2}) \text{ when } ||u||_{L^2} \to 0,
\end{equation}
\begin{equation} \label{approx_compk_O1}
\left\| \left(  
\langle \widetilde{\psi}(T) , \psi_j(T) \rangle 
- \omega_j \langle \mu \varphi_1 , \varphi_j \rangle \int_0^T s(t) e^{i \omega_j t} dt 
\right)_{j \in \JJ-\{1\}}  \right\|_{h^1}
= o(||s||_{L^2}) \text{ when } ||u||_{L^2} \to 0.
\end{equation}
\end{Prop}

The proof will use the following lemma which is a straightforward adaptation of \cite[Lemma 1]{KB-CL}. Its proof is postponed to Appendix \ref{appendix_lemme_KB-CL}.
\begin{Lem} \label{lemme:KB-CL}
Let $T>0$ and $f \in L^2((0,T),H^1)$. The function $F(t) := \int_0^t e^{iA \tau} f(\tau) d \tau$ belongs to
$C^0([0,T],H^1_0)$ and satisfies
\begin{equation*}
||F||_{L^{\infty}((0,T),H^1_0)} \leqslant c_1(T) ||f||_{L^2((0,T),H^1)}
\end{equation*}
where $c_1(T)>0$.
\end{Lem}

\noindent \textbf{Proof of Proposition \ref{Prop:approx_Comp1k}:}
Let $T>0$. We work with functions $u \in L^2((0,T),\R)$ such that $||u||_{L^2}<1$.
\\

\textit{First step : We prove that $||\widetilde{\psi}-\psi_1||_{L^{\infty}((0,T),H^1_0)} = O(||s||_{L^2})$ when $||u||_{L^2} \to 0$.}

\noindent
From Proposition \ref{WP-CYpb}, we know that $\psi \in C^0([0,T],H^3_{(0)})$ and
\begin{equation} \label{estimation_psi}
||\psi||_{L^{\infty}((0,T),H^3_{(0)})} \leqslant C.
\end{equation}
We deduce from (\ref{def:psitilde}) that $\widetilde{\psi} \in C^0([0,T],H^3 \cap H^1_0)$ and
\begin{equation} \label{estimation_psitilde}
||\widetilde{\psi}||_{L^{\infty}((0,T),H^3 \cap H^1_0)} \leqslant \tilde{C}.
\end{equation}
By Lemma \ref{lemme:KB-CL} the following equality holds in $H^1_0(0,1)$, for every $t \in [0,T]$
\begin{equation} \label{eq:psitilde}
\widetilde{\psi}(t) = \psi_1(t) - \int_0^t e^{iA(t-\tau)} \big[ s(\tau)( 2 \mu' \partial_x \widetilde{\psi}(\tau) + \mu'' \widetilde{\psi}(\tau) ) + i s(\tau)^2 \mu'^2 \widetilde{\psi}(\tau) \big] d \tau,
\end{equation}
and
\begin{align*}
|| \widetilde{\psi}- \psi_1||_{L^{\infty}((0,T),H^1_0)} 
\leqslant C(T) \Big[& ||s||_{L^2(0,T)} ||2 \mu' \partial_x \widetilde{\psi} + \mu'' \widetilde{\psi}||_{L^{\infty}((0,T),H^1)} 
\\
&+ ||s||_{L^2(0,T)}^2 ||\mu'^2 \widetilde{\psi}||_{L^{\infty}((0,T),H^1_0)} \Big].
\end{align*}
This inequality, together with (\ref{estimation_psitilde}) ends the first step.
\\

\textit{Second step : We prove that $||\widetilde{\psi}-\psi_1-\widetilde{\Psi}||_{L^{\infty}((0,T),H^1_0)} = o(||s||_{L^2})$ when $||u||_{L^2} \to 0$}

\noindent
where $\widetilde{\Psi}(t,x)$ is defined by
\begin{equation} \label{def:Psitilde}
\Psi(t,x) = \widetilde{\Psi}(t,x) + i s(t) \mu(x) \psi_1(t,x)
\end{equation}
and $\Psi$ is the solution of (\ref{Schro_lin_eq}). From Proposition \ref{WP-CYpb} (applied to system (\ref{Schro_lin_eq})), we know that $\Psi \in C^0([0,T],H^3_{(0)})$. We deduce from (\ref{def:Psitilde}) that
$\widetilde{\Psi} \in C^0([0,T],H^3 \cap H^1_0)$. Note that $\widetilde{\Psi}$ is a weak solution of
\begin{equation} \label{eq:Psitilde}
\left\lbrace
\begin{aligned}
&i \partial_t \widetilde{\Psi} = -\partial^2_x \widetilde{\Psi} - is(t)\big[ 2 \mu' \partial_x \psi_1 + \mu'' \psi_1\big],
\\
&\widetilde{\Psi}(t,0) = \widetilde{\Psi}(t,1) =0,
\\
&\widetilde{\Psi}(0,x) = 0.
\end{aligned}
\right.
\end{equation}
By Lemma \ref{lemme:KB-CL}, the following equality holds in $H^1_0(0,1)$, for every $t \in[0,T]$
\begin{equation} \label{Duhamel_Psitilde}
\widetilde{\Psi}(t) = - \int_0^t e^{iA(t-\tau)} s(\tau) \Big[ 2 \mu' \partial_x \psi_1(\tau) + \mu'' \psi_1(\tau) \Big] d \tau.
\end{equation}
Subtracting this relation to (\ref{eq:psitilde}) and applying Lemma \ref{lemme:KB-CL}, we get
\begin{align*}
||\widetilde{\psi}-\psi_1-\widetilde{\Psi}||_{L^{\infty}((0,T),H^1_0)} \leqslant & \,
C(T) \big( ||s||_{L^2(0,T)} || 2 \mu' \partial_x(\psi_1 -\widetilde{\psi}) + \mu'' (\psi_1 - \widetilde{\psi})||_{L^{\infty}((0,T),H^1)}
\\
&+ ||s||_{L^2(0,T)}^2 || \mu'^2 \widetilde{\psi}||_{L^{\infty}((0,T),H^1_0)} \big).
\end{align*}
We deduce from (\ref{estimation_psitilde}) the existence of a constant $C >0$ (independent of $u$) such that
\begin{equation*}
||\widetilde{\psi}-\psi_1-\widetilde{\Psi}||_{L^{\infty}((0,T),H^1_0)} \leqslant C \big( ||s||_{L^2} ||\widetilde{\psi}-\psi_1||_{L^{\infty}((0,T),H^2)} + ||s||_{L^2}^2 \big).
\end{equation*}
Thus, to end the proof of the second step, we only need to prove that
\begin{equation} \label{cv}
||\widetilde{\psi}-\psi_1||_{L^{\infty}((0,T),H^2)} \to 0 \; \text{ when } ||u||_{L^2} \to 0.
\end{equation}
Using (\ref{def:psitilde}) and  (\ref{estimation_psi}), we get
\begin{align*}
||\widetilde{\psi}-\psi_1||_{L^{\infty}((0,T),H^2)} 
&\leqslant ||(e^{is(t)\mu}-1) \psi||_{L^{\infty}((0,T),H^2)} + ||\psi - \psi_1||_{L^{\infty}((0,T),H^2)}
\\
&\leqslant C ||s||_{L^{\infty}(0,T)} + ||\psi-\psi_1||_{L^{\infty}((0,T),H^2)}.
\end{align*}
Thus (\ref{cv}) is a consequence of Proposition \ref{WP-CYpb} (applied to (\ref{Schro:syst})(\ref{IC})).
\\

\textit{Third step : We prove that $||\widetilde{\psi}-\psi_1-\widetilde{\Psi}-\widetilde{\xi}||_{L^{\infty}((0,T),L^2)} = o(||s||_{L^2}^2)$ when $||u||_{L^2} \to 0$}

\noindent
where $\widetilde{\xi}(t,x)$ is defined by
\begin{equation} \label{def:xitilde}
\xi(t,x) =\widetilde{\xi}(t,x) + is(t) \mu(x) \widetilde{\Psi}(t,x) - s(t)^2 \mu(x)^2 \psi_1(t,x)
\end{equation}
and $\xi$ is the solution of (\ref{Schro_O2_eq}). Note that $\widetilde{\xi}$ is a weak solution of
\begin{equation} \label{eq:xitilde}
\left\lbrace
\begin{aligned}
&i \partial_t \widetilde{\xi} = - \partial^2_x \widetilde{\xi} -  i s(t) \big[ 2 \mu'(x) \partial_x \widetilde{\Psi} + \mu''(x) \widetilde{\Psi} \big] + s(t)^2 \mu'(x)^2 \psi_1,
\\
&\widetilde{\xi}(t,0) = \widetilde{\xi}(t,1) =0,
\\
&\widetilde{\xi}(0,x) = 0.
\end{aligned}
\right.
\end{equation}
Thus, the following equation holds in $L^2(0,1)$ for every $t\in[0,T]$
\begin{equation} \label{Duhamel_xitilde}
\widetilde{\xi}(t) = -\int_0^t e^{iA(t-\tau)} \Big[ s(\tau) (2\mu' \partial_x \widetilde{\Psi}(\tau) + \mu'' \widetilde{\Psi}(\tau)) + is(\tau)^2 \mu'^2 \psi_1(\tau) \Big] d \tau.
\end{equation}
Using (\ref{eq:psitilde}) and (\ref{Duhamel_Psitilde}) we deduce that
\begin{align*}
(\widetilde{\psi}-\psi_1-\widetilde{\Psi}-\widetilde{\xi})(t) 
=& -\int_0^t e^{iA(t-\tau)} \Big[ s(\tau) \Big( 2\mu' \partial_x (\widetilde{\psi}-\psi_1-\widetilde{\Psi})(\tau) + \mu''(\widetilde{\psi}-\psi_1-\widetilde{\Psi})(\tau) \Big) 
\\
& + is(\tau)^2 \mu'^2 (\widetilde{\psi}-\psi_1)(\tau) \Big] d \tau
\end{align*}
in $L^2(0,1)$ for every $t \in [0,T]$. Thus,
\begin{equation*}
||(\widetilde{\psi}-\psi_1-\widetilde{\Psi}-\widetilde{\xi})(t)||_{L^2} \leqslant C \int_0^t |s(\tau)| \, ||(\widetilde{\psi}-\psi_1-\widetilde{\Psi})(\tau)||_{H^1} + |s(\tau)|^2 || (\widetilde{\psi}-\psi_1)(\tau) ||_{L^2} d \tau
\end{equation*}
Taking into account the first and second step, we get the conclusion of the third step.
\\

\textit{Fourth step : Proof of (\ref{approx_comp1_O2}).} We deduce from (\ref{Duhamel_Psitilde}) and (\ref{Duhamel_xitilde}) that
\begin{equation*}
\Im \langle \widetilde{\Psi}(T) , \varphi_K e^{-i \lambda_1 T} \rangle =0,
\quad
\Im \langle \widetilde{\xi}(T) , \varphi_K e^{-i \lambda_1 T} \rangle = \QQ_{K,T}(s).
\end{equation*}
Using the third step, we get
\begin{align*}
\Big|  \Im \langle \widetilde{\psi}(T) , \varphi_K e^{-i \lambda_1 T} \rangle - \QQ_{K,T}(s) \Big|
&= \Big| \Im \langle (\widetilde{\psi}-\psi_1-\widetilde{\Psi}-\widetilde{\xi})(T), \varphi_K e^{-i \lambda_1 T} \rangle \Big|
\\
& \leqslant  || (\widetilde{\psi}-\psi_1-\widetilde{\Psi}-\widetilde{\xi})(T) ||_{L^2} 
\\
& = o(||s||_{L^2}^2) \; \text{ when } ||u||_{L^2} \to 0.
\end{align*}
\\

\textit{Fifth step : Proof of (\ref{approx_O1}).} We deduce from (\ref{Duhamel_Psitilde}) and the relation $\langle 2 \mu' \varphi_1' + \mu'' \varphi_1, \varphi_1 \rangle = 0$ that $\Im \langle \widetilde{\Psi}(T),\psi_1(T) \rangle =0$. Thus, the second step gives 
\begin{align*}
\Big| \Im \langle \widetilde{\psi}(T), \psi_1(T) \rangle \Big| &= \Big| \Im \langle (\widetilde{\psi}-\psi_1-\widetilde{\Psi})(T), \psi_1(T) \rangle \Big|
\\
&= o(||s||_{L^2}) \; \text{ when } ||u||_{L^2} \to 0.
\end{align*}
\\

\textit{Sixth step : Proof of (\ref{approx_compk_O1}).} We deduce from (\ref{Duhamel_Psitilde}) that
\begin{equation*}
\langle \widetilde{\Psi}(T), \psi_j(T) \rangle = \omega_j \langle \mu \varphi_1, \varphi_j \rangle \int_0^T s(t) e^{i \omega_j t} d t, \quad  \forall j \in \N^*-\{1\}.
\end{equation*}
Using the second step, we get
\begin{align*}
& \Big| \Big| \Big( \langle \widetilde{\psi}(T),\psi_j(T) \rangle - \omega_j \langle \mu \varphi_1, \varphi_j \rangle \int_0^T s(t) e^{i \omega_j t} d t \Big)_{j \in \JJ-\{1\} } \Big| \Big|_{h^1}
\\
&= \Big| \Big| \big( \langle (\widetilde{\psi}- \psi_1 - \widetilde{\Psi})(T), \psi_j(T) \big)_{j \in \JJ-\{1\} } \Big| \Big|_{h^1}
\\
&\leqslant C || (\widetilde{\psi}- \psi_1 - \widetilde{\Psi})(T) ||_{H^1_0}
\\
&= o(||s||_{L^2}) \; \text{ when } ||u||_{L^2} \to 0.
\end{align*}
This ends the proof of Proposition \ref{Prop:approx_Comp1k}.
\hfill $\blacksquare$
\\

\noindent \textbf{Proof of Theorem \ref{Thm:aux}:}
One may assume that $A_K>0$, $\alpha_K=1$. Let $T<T_K^*$.
Working by contradiction, we assume that, for every $\epsilon>0$, 
there exists $s_\epsilon \in H^1((0,T),\mathbb{R})$ with $s_\epsilon(0)=0$ and $\|\dot{s}_\epsilon\|_{L^2}<\epsilon$ 
such that the solution $\widetilde{\psi}_{\epsilon}$ of  (\ref{Schro:syst_aux}) satisfies 
\begin{equation} \label{target_imp}
\widetilde{\psi}_\epsilon(T,.) = ( \sqrt{1-\delta_\epsilon^2} \varphi_1 + i \delta_\epsilon \varphi_K ) 
e^{i \theta_\epsilon \mu(.)} e^{-i\lambda_1 T}
\end{equation}
for some $\delta_\epsilon>0$ and $\theta_\epsilon \in \mathbb{R}$.
Then $\theta_\epsilon, \delta_\epsilon \rightarrow 0$ when $\epsilon \rightarrow 0$.
\\

\noindent \emph{First step: We prove that
\begin{equation} \label{Equiv_th_de}
|\theta_\epsilon|+|\delta_\epsilon|=\underset{\epsilon \rightarrow 0}{O} ( \|s_\epsilon\|_{L^2} ).
\end{equation}}
Using (\ref{target_imp}) and the assumption $\langle \mu \varphi_1, \varphi_K \rangle =0$, we have
\begin{align*}
& \frac{1}{2} \| (\widetilde{\psi}_\epsilon-\psi_1)(T)\|_{L^2(0,1)}^2 \\
= & 1 - \Re \int_0^1 \widetilde{\psi}_\epsilon(T,x) \overline{\psi_1(T,x)} dx 
\\
= & 1 - \int_0^1 \Big( 
\sqrt{1-\delta_\epsilon^2} \varphi_1(x)^2 \cos[\theta_\epsilon \mu(x)]    
-\delta_\epsilon \varphi_1(x)\varphi_K(x)\sin[\theta_\epsilon \mu(x)] 
\Big)dx
\\
= & 1 - 
\left( 1-\frac{\delta_\epsilon^2}{2} + O(\delta_\epsilon^4) \right) 
\left( 1-\frac{\theta_\epsilon^2}{2}\|\mu\varphi_1\|^2+O(\theta_\epsilon^4) \right)
+ \underset{\epsilon \rightarrow 0}{O}(\delta_\epsilon \theta_\epsilon^3)
\\ = &
\frac{\delta_\epsilon^2}{2} + \frac{\theta_\epsilon^2}{2}\|\mu\varphi_1\|^2 +
\underset{\epsilon \rightarrow 0}{O}( \delta_\epsilon^4 + \theta_\epsilon^4 + \delta_\epsilon \theta_\epsilon^3 ).
\end{align*}
As proved, in Proposition \ref{Prop:approx_Comp1k},
\begin{equation*}
||\widetilde{\psi}-\psi_1||_{L^{\infty}((0,T),H^1_0)} = O(||s||_{L^2}) \; \text{ when } ||u||_{L^2} \to 0.
\end{equation*}
This concludes the first step.
\\

\noindent \emph{Second step: We prove that
\begin{equation} \label{step2_aux}
\Im \langle  \widetilde{\psi}_\epsilon(T) , \varphi_K e^{-i\lambda_1 T} \rangle
=\delta_\epsilon + \underset{\epsilon \rightarrow 0}{o} ( \|s_\epsilon\|_{L^2}^2 ).
\end{equation}}
Using (\ref{target_imp}) and the assumption $\langle\mu\varphi_1,\varphi_K\rangle=0$, we get
\begin{align*}
\Im \langle  \widetilde{\psi}_\epsilon(T) , \varphi_K e^{-i\lambda_1 T} \rangle
& = 
\int_0^1 \Big(
\sqrt{1-\delta_\epsilon^2} \varphi_1(x) \varphi_K(x) \sin[\theta_\epsilon \mu(x)] 
+\delta_\epsilon  \varphi_K(x)^2 \cos[\theta_\epsilon \mu(x)] 
\Big) dx
\\ 
& = 
\left( 1 + \underset{\epsilon\rightarrow 0}{O}(\delta_\epsilon^2) \right) \underset{\epsilon\rightarrow 0}{O}(\theta_\epsilon^3)
+ \delta_\epsilon \left( 1  + \underset{\epsilon\rightarrow 0}{O}(\theta_\epsilon^2) \right)
\\
& =
\delta_\epsilon + \underset{\epsilon\rightarrow 0}{O}( \theta_\epsilon^3 + \delta_\epsilon \theta_\epsilon^2 ).
\end{align*}
Thus, (\ref{step2_aux}) is a consequence of (\ref{Equiv_th_de}).
\\

\noindent \emph{Third step: Conclusion.} Using (\ref{step2_aux}), (\ref{approx_comp1_O2}) and (\ref{QKT_coercive}), we get
\begin{align*}
0 < \delta_\epsilon 
   & = \Im \langle  \widetilde{\psi}_\epsilon(T) , \varphi_K e^{-i\lambda_1 T} \rangle 
+ \underset{\epsilon\rightarrow 0}{o}(  \|s_\epsilon\|_{L^2}^2  )
\\ & = \mathcal{Q}_{K,T}(s_\epsilon) + \underset{\epsilon\rightarrow 0}{o}( \|s_\epsilon\|_{L^2}^2 )
\\ & \leqslant - \frac{A_K}{4} \|s_\epsilon\|_{L^2}^2 + \underset{\epsilon\rightarrow 0}{o}( \|s_\epsilon\|_{L^2}^2),
\end{align*}
which is impossible when $\epsilon$ is small. \hfill  $\blacksquare$

\section{Local controllability in large time}
\label{sec:thm:control_Tlarge}

The goal of this section is to prove Theorem \ref{thm:control_Tlarge}.

\subsection{Preliminary}

The goal of this section is the proof of the following result.

\begin{Prop} \label{Prop:mu'pm}
Let $\mu \in H^3((0,1),\mathbb{R})$ be such that $\mu'(0)\pm\mu'(1)\neq 0$.
\begin{enumerate}
\item Then, $N:=\sharp\{k \in \mathbb{N}^* ; \langle \mu \varphi_1,\varphi_k \rangle=0\}$ is finite.
\item Let $K_1<...<K_N \in \mathbb{N}^*$ be such that $\langle\mu\varphi_1,\varphi_{K_j}\rangle=0$ for $j=1,...,N$.
Then, for every $j \in \{1,...,N\}$ and $T>0$ there exists $v \in V_T$ such that $Q_{K_j,T}^2(v) \neq 0$.
\item There exists $c>0$ such that
\begin{equation} \label{hyp_mu_2.2}
|\langle \mu \varphi_1,\varphi_k\rangle | \geqslant \frac{c}{k^3}, 
\forall k \in \mathbb{N}^*-\{K_1,...,K_N\}.
\end{equation}	
\end{enumerate}
\end{Prop}

We recall that $Q_{K_j,T}^2$ is defined in (\ref{def:Q2KT})(\ref{def:h2K}), and $V_T$ in (\ref{def:VT}).
For the proof of Proposition \ref{Prop:mu'pm}, we need the following preliminary result.

\begin{Prop} \label{Prop:Q2<>0}
Let $\mu\in H^3((0,1),\mathbb{R})$ and $K\in \mathbb{N}^*$ be such that 
$\langle \mu \varphi_K,\varphi_n\rangle\langle\mu\varphi_n,\varphi_1\rangle \neq 0$
for some $n \in \mathbb{N}^*$. The following statements are equivalent.
\begin{itemize}
\item There exists $T^*>0$ such that, for every $T<T^*$, $Q_{K,T}^2 \equiv 0$ on $V_T$.
\item The support of the sequence $(\langle \mu \varphi_K,\varphi_j\rangle\langle\mu\varphi_j,\varphi_1\rangle)_{j \in \mathbb{N}^*}$
is contained in the finite set 
$$\{ j_* \in \JJ \cap [1,K] ; \exists k_* \in \JJ \cap [1,K] ,\lambda_{j_*}-\lambda_1=\lambda_K-\lambda_{k_*} \},$$
and for every $j_*,k_* \in \JJ \cap [1,K]$ such that $\lambda_{j_*}-\lambda_1=\lambda_K-\lambda_{k_*}$
then 
$\langle \mu \varphi_K,\varphi_{j_*}\rangle\langle\mu\varphi_{j_*},\varphi_1\rangle
=\langle \mu \varphi_K,\varphi_{k_*}\rangle\langle\mu\varphi_{k_*},\varphi_1\rangle$.
\end{itemize}
\end{Prop}

\noindent \textbf{Proof of Proposition \ref{Prop:Q2<>0}:}
To simplify the notation of this proof, we write $Q_T$ and $h$, instead of $Q^2_{K,T}$ and $h^2_K$.
Let us assume that $Q_{T} \equiv 0$ on $V_T$, for every $T<T^*$.
Then $\nabla Q_T (v) \perp V_T$, for every $v \in V_T$ and $T<T^*$.
Easy computations show that, for $v \in V_T$,
\begin{equation} \label{nablaQT=}
\nabla Q_{T}(v):t \mapsto \int_0^t v(\tau)h(t,\tau)d\tau+\int_t^{T}v(\tau)h(\tau,t)d\tau=
\int_t^T v(\tau)[h(\tau,t)-h(t,\tau)]d\tau.
\end{equation}

\noindent \emph{First step: We prove that $\nabla Q_{T}(v)=0, \forall v \in V_T$.}
Let $T$, $T_1$ be such that $0<T<T_1<T^*$ and $v \in V_{T_1}$ supported on $(0,T)$.
Since $\nabla Q_{T_1}(v) \perp V_{T_1}$, there exists a unique sequence $(\alpha_k)_{k \in \mathbb{Z}-\{0\}} \in l^2$ such that
$$\nabla Q_{T_1}(v) = \sum_{k=1}^{\infty}  \alpha_k e^{i(\lambda_k-\lambda_1)t} 
+ \sum_{k=2}^{\infty} \alpha_{-k}  e^{-i(\lambda_k-\lambda_1)t} 
\quad \text{ in } L^2((0,T_1),\mathbb{R})$$
(decomposition on a Riesz-basis).
We have $\nabla Q_{T_1}(v) \equiv 0$ on $(T,T_1)$ because $v$ is supported on $(0,T)$ (see(\ref{nablaQT=})).
Using Ingham inequality on $(T,T_1)$ we get $\alpha_k=0,\forall k$ (see Proposition \ref{Cor:haraux1} in Appendix).
\\

\noindent \emph{Second step: We prove that $V_{T}|_{(t,T)}=L^2(t,T)$.}
Let $T \in (0,T^*)$ and $t \in (0,T)$ be fixed. Let $v \in L^2((t,T),\R)$. We define $d_j := 0$, for $j \in \N^*- \JJ$ and
\begin{equation*}
d_j := - \int_t^T v(\tau) e^{i\omega_j \tau} d \tau, \quad \text{for } j \in \JJ.
\end{equation*}
Thus, $d=(d_j)_{j\in \N^*} \in \ell^2_r(\N^*,\C)$ and Proposition \ref{Cor:haraux1} imply that there exists $\tilde{v} \in L^2((0,t),\R)$ such that
\begin{equation*}
\int_0^t \tilde{v}(\tau) e^{i \omega_j \tau} d\tau = d_j = - \int_t^T v(\tau) e^{i\omega_j \tau} d \tau, \quad \forall j \in \JJ.
\end{equation*}
Then if we extend $\tilde{v}$ on $(t,T)$ by setting $\tilde{v}_{|(t,T)} =v$, it comes that $\tilde{v} \in V_T$. Thus, $V_{T}|_{(t,T)}=L^2(t,T)$.
\\

\noindent \emph{Third step: We prove that $h(\tau,t)=h(t,\tau), \forall t,\tau \in [0,T^*]$.} 
Using the first step, we get
\begin{equation} \label{nablaQ=0_bis}
\int_t^T v(\tau)[h(\tau,t)-h(t,\tau)]d\tau = 0, \quad \forall 0<t<T<T^*, \forall v \in V_T.
\end{equation}
Using the second step, we deduce from (\ref{nablaQ=0_bis}) that
$\tau \mapsto h(t,\tau)-h(\tau,t)$ vanishes in $L^2(t,T)$,
for every $0<t<T<T^*$. This gives the conclusion because the function $(t,\tau) \mapsto h(\tau,t)-h(t,\tau)$ is continuous.
\\

\noindent \emph{Fourth step: Conclusion.}
Let $k^* \in \mathbb{N}^*$ be such that
$b_{k^*}:=\langle\mu\varphi_K,\varphi_{k^*}\rangle\langle\mu\varphi_{k^*},\varphi_1\rangle \neq 0$.
The equality $h(t,\tau)-h(\tau,t)=0$ with $\tau=0$ gives
\begin{equation} \label{h(t,0)=h(0,t)}
b_{k^*} e^{i(\lambda_K-\lambda_{k^*})t}
= \sum\limits_{j \in \JJ} b_j e^{i(\lambda_j-\lambda_1)t} - \sum\limits_{k \in \JJ-\{k^*\}} b_k e^{i(\lambda_K-\lambda_k)t}.
\end{equation}
The equality $\frac{d}{d\tau}[h(t,\tau)-h(\tau,t)]=0$ with $\tau=0$ gives
\begin{equation} \label{der_h(t,0)=der_h(0,t)}
(\lambda_{k^*}-\lambda_1) b_{k^*}  e^{i(\lambda_K-\lambda_{k^*})t}
= \sum\limits_{j \in \JJ} (\lambda_K-\lambda_j) b_j  e^{i(\lambda_j-\lambda_1)t} 
- \sum\limits_{k \in \JJ-\{k^*\}} (\lambda_k-\lambda_1) b_k e^{i(\lambda_K-\lambda_k)t}.
\end{equation}
Thus, an obvious linear combination of (\ref{h(t,0)=h(0,t)}) and (\ref{der_h(t,0)=der_h(0,t)}) leads to
$$ 0= \sum\limits_{j \in \JJ} \big( (\lambda_K-\lambda_j) - (\lambda_{k^*}-\lambda_1) \big) b_j  e^{i(\lambda_j-\lambda_1)t} 
- \sum\limits_{k \in \JJ-\{k^*\}} \big( (\lambda_k-\lambda_1) - (\lambda_{k^*}-\lambda_1) \big) b_k e^{i(\lambda_K-\lambda_k)t}.$$
In the right hand side of the previous equality,
the frequencies $(\lambda_j-\lambda_1)$ are $\geqslant 0$ for every $j \in \JJ$,
while the frequencies $(\lambda_K-\lambda_k)$ are negative for every $k>K$.
Thus, for every $k>K$ the frequency $(\lambda_K-\lambda_k)$ appears only one time
in the right hand side of the previous equality. 
The uniqueness of the decomposition on a Riesz basis gives
$$(\lambda_{k^*}-\lambda_1) b_k = (\lambda_k-\lambda_1) b_k, \forall k>K \text{ with } k \neq k^*.$$
Thus, $b_k=0$, $\forall k \in \JJ-\{k^*\}$ with $k>K$.
Coming back to (\ref{h(t,0)=h(0,t)}), we only have a finite sum in the right hand side,
over $j \in \JJ$ with $j \leqslant K$ and over $k \in \JJ-\{k^*\}$ with $k \leqslant K$.
We deduce the existence of a unique $j^* \in \JJ$ with $j^* \leqslant K$ such that
$\lambda_K-\lambda_{k^*}=\lambda_{j^*}-\lambda_1$ and $b_{k^*}=b_{j^*}$.

Reciprocally, let 
$\alpha:=\lambda_K-\lambda_{k^*}=\lambda_{j^*}-\lambda_1$,
$\beta:=\lambda_K-\lambda_{j^*}=\lambda_{k^*}-\lambda_1$.
Then $h(t,\tau):=b_{k^*}[e^{i[\alpha t +\beta \tau]}+e^{i[\beta t + \alpha \tau]} ]$,
satisfies $h(t,\tau)=h(\tau,t)$ and $\nabla Q_{T} \equiv 0$ on $V_T$, for every $T>0$. 
By linearity, the same conclusion holds when $h$ is a finite sum of such terms. 

\hfill  $\blacksquare$

\noindent \textbf{Proof of Proposition \ref{Prop:mu'pm}:} 
Performing three integrations by part and using the Riemann-Lebesgue Lemma, we get for every $K$ and $n$ in $\N^*$,
\begin{equation} \label{3IPP}
\langle\mu\varphi_K,\varphi_n\rangle
= \frac{4K[(-1)^{K+n}\mu'(1)-\mu'(0)]}{n^3 \pi^2} + \underset{n \rightarrow + \infty}{o}\left( \frac{1}{n^3} \right).
\end{equation}
Thus, for $n$ large enough $\langle \mu \varphi_1,\varphi_n\rangle \neq 0$. 
This proves the first and third statements of Proposition \ref{Prop:mu'pm}.

Let $j \in \{1,...,N\}$. Using (\ref{3IPP}), we have simultaneously
$\langle \mu \varphi_1,\varphi_n\rangle \neq 0$ and 
$\langle \mu \varphi_{K_j},\varphi_n\rangle \neq 0$ for arbitrarily large values of $n$.
Thus, Proposition \ref{Prop:Q2<>0} gives the conclusion. 

\hfill  $\blacksquare$

\subsection{Strategy for the proof of Theorem \ref{thm:control_Tlarge}}

Until the end of Section \ref{sec:thm:control_Tlarge}, we fix 
$\mu \in H^3((0,1),\mathbb{R})$ such that $\mu'(1)\pm\mu'(0)\neq 0$,
$N \in \mathbb{N}$ and  $K_1,...,K_N \in \mathbb{N}^*$ as in Proposition \ref{Prop:mu'pm}.
To simplify the notations, we assume that $K_1=1$. We define the space
\begin{equation} \label{def:H}
\mathcal{H} := \text{Span}_{\C} \Big( \psi_k(T) , k \in \mathbb{N}^*-\{K_1,...,K_N\} \Big),
\end{equation}
and, for $j=1,...,N$ the space
\begin{equation} \label{def:Mj}
M^j := \left\lbrace \begin{array}{l}
\text{Span}_{\C} \left( \psi_{K_j}(T) \right) \text{ if } K_j \neq 1,\\
i \text{Span}_{\R} (\psi_1(T)) \text{ if } K_j=1.
\end{array}\right.
\end{equation}
Let
\begin{equation} \label{def:M}
M := \bigoplus_{j=1}^N M^j. 
\end{equation}
The global strategy relies on power series expansion 
of the solutions to the second order as in \cite{Crepeau-Cerpa} (see also \cite{JMC-book}). 
In Section \ref{subsec:H}, we prove the local exact controllability 'in $\mathcal{H}$', with a first order strategy.
Then, in Section \ref{subsec:O2}, we prove that any direction in $M$ is reached with the second order term.
Finally, in Section \ref{subsec:Fix_point}, we conclude with a fixed point argument.

\subsection{Controllability in $\mathcal{H}$ in arbitrarily small time}
\label{subsec:H}

We introduce the orthogonal projection
\begin{equation} \label{def:P1T}
\begin{array}{|cccl}
\mathcal{P}_{T}:   & L^2(0,1) & \rightarrow & \mathcal{H} \\
                   &   \psi    & \mapsto     & \psi - \sum\limits_{j=1}^N \langle \psi , \psi_{K_j}(T) \rangle \psi_{K_j}(T)
\end{array}
\end{equation}
The goal of this section is the proof of the following result.

\begin{thm} \label{thm:control_H}
Let $T_1,T>0$ be such that $T_1<T$. There exists $\delta_1>0$ and a $C^1$-map
$\Gamma_{[T_1,T]}: \Omega_{T_1} \times \Omega_T  \rightarrow  L^2((T_1,T),\mathbb{R})$
where
$$\Omega_{T_1} := \{ \psi_0 \in \mathcal{S} \cap H^3_{(0)}(0,1) ; \| \psi_0-\psi_1(T_1) \|_{H^3_{(0)}} < \delta_1 \},$$
$$\Omega_T:= \{ \widetilde{\psi}_f \in \mathcal{H} \cap H^3_{(0)}(0,1) ; \| \widetilde{\psi}_f - \mathcal{P}_T[\psi_1(T)] \|_{H^3_{(0)}} < \delta_1 \}$$
such that $\Gamma_{[T_1,T]} (\psi_1(T_1),\mathcal{P}_T[\psi_1(T)])=0$ and for every 
$(\psi_0, \widetilde{\psi}_f) \in \Omega_{T_1} \times \Omega_T$,
the solution of (\ref{Schro:syst}) with 
initial condition $\psi(T_1)=\psi_0$ and control $u:=\Gamma_{[T_1,T]} (\psi_0,\widetilde{\psi}_f)$ satisfies 
$\mathcal{P}_{T}[\psi(T)]=\widetilde{\psi}_f$.
\end{thm}

This theorem may be proved exactly as Theorem \ref{Thm:KB-CL} in \cite{KB-CL}.
We recall the main steps of the proof because several intermediate results will also be used in the end of this article. 
To simplify the notations, we take $T_1=0$.
\\

By Proposition \ref{WP-CYpb}, we can consider the map
\begin{equation} \label{def:Theta1T}
\begin{array}{|cccc}
\Theta_{T}: & [\mathcal{S} \cap H^3_{(0)}(0,1)]  \times L^2((0,T),\mathbb{R}) 
& \rightarrow & [\mathcal{S} \cap H^3_{(0)}(0,1)]  \times [\mathcal{H} \cap H^3_{(0)}(0,1)] 
\\
& (\psi_0,u) & \mapsto & (\psi_0,\mathcal{P}_{T}[\psi(T)])
\end{array}
\end{equation}
where $\psi$ is the solution of (\ref{Schro:syst})(\ref{IC_qlq}).
Then Theorem \ref{thm:control_H} corresponds to the local surjectivity of the nonlinear map $\Theta_{T}$
around the point $(\varphi_1,0)$, that will be proved thanks to the inverse mapping theorem.
Thus, the first property required is the $C^1$-regularity of $\Theta_{T}$, 
which is a consequence of \cite[Proposition 3]{KB-CL}.

\begin{Prop} \label{Prop:C1}
Let $T>0$ and $\mu \in H^{3}((0,1),\mathbb{R})$.
The map $\Theta_{T}$ defined by (\ref{def:Theta1T}) is $C^1$.
Moreover, for every 
$\psi_0, \Psi_0 \in H^3_{(0)}(0,1)$, 
$u, v \in L^2((0,T),\mathbb{R})$, we have
\begin{equation} \label{dTheta1T}
d\Theta_{T}(\psi_0,u).(\Psi_0,v)=(\Psi_0,P_{T}[\Psi(T)])
\end{equation}
where $\Psi$ is the weak solution of the linearized system
\begin{equation} \label{proofC1:L}
\left\lbrace\begin{array}{ll}
i \partial_t \Psi = - \partial_x^2 \Psi - u(t) \mu(x) \Psi - v(t) \mu(x) \psi,&  (t,x) \in (0,T)\times(0,1),\\
\Psi(t,0)=\Psi(t,1)=0,                                                              & t \in (0,T), \\
\Psi(0,x)=\Psi_0,                                                                   & x \in (0,1)
\end{array}\right.
\end{equation} 
and $\psi$ is the solution of (\ref{Schro:syst})(\ref{IC_qlq}).
\end{Prop}

The second property required for the application of the inverse mapping theorem is the
the existence of a continuous right inverse for $d\Theta_{T}(\varphi_1,0)$, 
that may be proved exactly as \cite[Proposition 4]{KB-CL}
(it is a consequence of Proposition \ref{Cor:haraux1} in Appendix).

\begin{Prop} \label{Cont-Lin-H3L2}
Let $T>0$ and $\mu \in H^{3}((0,1),\mathbb{R})$ be such that (\ref{hyp_mu_2.2}) holds. 
The linear map 
$$d\Theta_{T}(\varphi_1,0):
[T_{\mathcal{S}} \varphi_1 \cap H^3_{(0)}]  \times L^2((0,T),\mathbb{R}) 
\rightarrow 
[T_{\mathcal{S}} \varphi_1 \cap H^3_{(0)}]  \times [\mathcal{H} \cap H^3_{(0)}]  $$
has a continuous right inverse
$$d\Theta_{T}(\varphi_1,0)^{-1}:
[T_{\mathcal{S}} \varphi_1 \cap H^3_{(0)}]  \times [\mathcal{H} \cap H^3_{(0)}] 
\rightarrow
[T_{\mathcal{S}} \varphi_1 \cap H^3_{(0)}]  \times L^2((0,T),\mathbb{R}).$$
\end{Prop}
Thus, Propositions \ref{Prop:C1} and \ref{Cont-Lin-H3L2} allow to apply the inverse mapping theorem to $\Theta_T$ 
at the point $(\varphi_1,0)$ and thus to prove Theorem \ref{thm:control_H}.

\subsection{Reaching the missed directions, at the second order, in large time.}
\label{subsec:O2}

The goal of this section is the proof of the following result.

\begin{Prop} \label{Prop:dir_perdue_co}
Let $T>T_\sharp$ where
\begin{equation} \label{def_Tmin}
T_\sharp := 
\left\lbrace \begin{aligned}
&2^{N-1} T_{min}^2 + \sum_{k=2}^N ((k-1) + 2^{k-2}) \frac{\pi}{\lambda_{K_k} - \lambda_1} & \text{ if } K_1 = 1, 
\\
&\sum_{k=1}^N  \frac{k\pi}{\lambda_{K_k} - \lambda_1}                                     & \text{ if } K_1 \neq 1. 
\end{aligned}\right.
\end{equation}
There exists a continuous map
$$\begin{array}{|cccl}
\Lambda_T: & M & \rightarrow & L^2((0,T),\mathbb{R})^2 \\
           & z & \mapsto     & (v,w)
\end{array}$$
such that, for every $z \in M$, the solutions $\Psi$ and $\xi$ of
(\ref{Schro_lin_eq}) and (\ref{Schro_O2_eq}) satisfy $\Psi(T)=0$ and $\xi(T)=z$.
\end{Prop}

In this statement, the quantity $T_{min}^2$ is defined as follows.

\begin{Lem} \label{Lem:Q1<>0}
The quantity 
$$T_{min}^2 := \inf\{ T>0 ; \exists v_{\pm} \in V_T \text{ such that } \widetilde{Q}_{1,T}^2(v_\pm)=\pm 1 \}$$
is well defined and belongs to $(0,2/\pi]$.
\end{Lem}

Let us recall that $K_1=1$ and $\widetilde{Q}_{1,T}^2$ and $V_T$ are defined in (\ref{def:Q2KTtilde}), (\ref{def:VT}).
\\

\noindent \textbf{Proof of Lemma \ref{Lem:Q1<>0}:} Let $T \geqslant 2/\pi$. Let $T_1^*$ be defined as in Lemma \ref{Lem:T_2}. If $v_- \in V_T-\{0\}$ is supported on $(0,T_1^*)$, then Lemma \ref{Lem:T_2} implies that $\widetilde{Q}_{1,T}^2(v)<0$.
Let $v_+(t):=\cos( \pi^2 t) 1_{[0,2/\pi]}(t)$. Then, explicit computations prove that $v_+ \in V_{T}$ and
$$\widetilde{Q}_{1,T}^2(v_+)=\sum\limits_{j=2}^\infty \langle \mu \varphi_1,\varphi_j \rangle^2 \frac{(j^2-1)}{\pi^3 j^2 (j^2-2)} >0. $$
\hfill $\blacksquare$

\subsubsection{Preliminary}

Our proof of Proposition \ref{Prop:dir_perdue_co} requires three preliminary results.
The first one consists in proving the existence of controls such that the projections 
of the second order term on the lost directions are non zero. 

\begin{Prop} \label{prop_dir_perdues}
Let $T>0$. For every $j\in\{1,...,N\}$, there exists $v_j, w_j \in L^2((0,T),\R)$ 
such that the associated solutions $\Psi^j$ and $\xi^j$ of (\ref{Schro_lin_eq}) and (\ref{Schro_O2_eq}) satisfy
\begin{align*}
&\Psi^j(T,\cdot) =0,
\\
&\lag \xi^j(T,\cdot) , \psi_{K_j}(T) \rag \neq 0, 
\\
&\lag \xi^j(T,\cdot) , \psi_k(T) \rag =0, \quad \forall k \in \mathbb{N}^*-\{K_1,...,K_N\}.
\tag{\theequation} \addtocounter{equation}{1}
\label{dir_perdue_projection}
\end{align*}
\end{Prop}

\noindent \textbf{Proof of Proposition \ref{prop_dir_perdues}:} Let $j \in \{1,...,N\}$.
By Proposition \ref{Prop:mu'pm}, there exists
$v_j \in V_T$ such that $Q^2_{K_j,T}(v_j) \neq 0$.
Using (\ref{CompK_O2}) we get $\lag \xi^j(T) , \psi_{K_j}(T) \rag = Q^2_{K_j,T}(v_j) \neq 0$.
As $v_j \in V_T$, (\ref{Psi_explicit}) and (\ref{def:VT}) imply $\Psi^j(T)=0$.
The equality (\ref{dir_perdue_projection}) is equivalent to the following trigonometric moment problem on $w_j$,
\begin{equation} \label{moment_pb_O2}
\int_0^T w_j(t) e^{i\omega_k t} dt = \frac{1}{\langle \mu \varphi_1 , \varphi_k \rangle}
\int_0^T v_j(t) \langle \mu \Psi^j(t),\varphi_k \rangle e^{i\lambda_k t} dt, 
\forall k \in \mathbb{N}^*-\{K_1,...,K_N\}.
\end{equation}
By (\ref{hyp_mu_2.2}) and \cite[Lemma 1]{KB-CL}, the right hand side 
belongs to $l^2$. Thus, Proposition \ref{Cor:haraux1} ensures the existence of a solution
$w_j \in L^2((0,T),\mathbb{R})$. \hfill  $\blacksquare$
\\

The second preliminary result for the proof of Proposition \ref{Prop:dir_perdue_co}
is a measure of the rotation of the null input solution, precised in the next statement.

\begin{Lem} \label{lemme_sans_controle}
Let $T, \tilde{T}, \theta>0$ be such that $0<T<T+\theta \leqslant \tilde{T}$,
$v,w \in L^2((0,T),\mathbb{R})$ and $v_\theta, w_\theta \in L^2((0,\tilde{T}),\mathbb{R})$ be defined by
\begin{equation*}
(v_\theta,w_\theta)(t) :=
\left\{\begin{array}{ll}
(0,0) & \quad \text{ if } t \in (0,\theta),
\\
(v,w)(t-\theta) & \quad  \text{ if } t \in (\theta,\theta+T),
\\
(0,0)  & \quad \text{ if } t \in (\theta+T,\tilde{T}).
\end{array}\right.
\end{equation*}
We denote by $(\Psi,\xi)$ and $(\Psi_\theta,\xi_\theta)$ the associated solutions 
of (\ref{Schro_lin_eq}) and (\ref{Schro_O2_eq}). Then, for every $k \in \N^*$
\begin{align*}
\lag \Psi_\theta(\tilde{T}) , \psi_k(\tilde{T}) \rag &= e^{i (\lambda_k - \lambda_1) \theta} \lag \Psi(T) , \psi_k(T) \rag,
\\
\lag \xi_\theta(\tilde{T}) , \psi_k(\tilde{T}) \rag &= e^{i (\lambda_k - \lambda_1) \theta} \lag \xi(T) , \psi_k(T) \rag.
\end{align*}
\end{Lem}

\begin{rk}
\label{rk_rotation}
Note that, for $k=1$, there is no rotation phenomenon.
\end{rk}

\noindent \textbf{Proof of Lemma \ref{lemme_sans_controle}:} We have
$$\Psi_\theta(t)= \left\lbrace \begin{array}{ll}
0                                          & \text{ for } 0<t<\theta,
\\
\Psi(t-\theta)e^{-i\lambda_1\theta}        & \text{ for } \theta<t<\theta+T,
\\
e^{-iA(t-\theta-T)} \Psi_\theta(\theta+T)  & \text{ for } \theta+T<t \leqslant \tilde{T},
\end{array}\right.$$
thus
$$\begin{array}{ll}
\Psi_\theta(\tilde{T}) 
& = \sum\limits_{k=1}^\infty \langle \Psi(T) , \varphi_k \rangle e^{-i\lambda_1 \theta} e^{-i\lambda_k(\tilde{T}-\theta-T)} \varphi_k
\\ 
& =\sum\limits_{k=1}^\infty \langle \Psi(T),\psi_k(T)\rangle e^{i(\lambda_k-\lambda_1)\theta} \psi_k(\tilde{T}).
\end{array}$$
The same relations hold for $\xi_\theta$.\hfill  $\blacksquare$
\\

The third preliminary result for the proof of Proposition \ref{Prop:dir_perdue_co}
is the non overlapping principle.

\begin{Prop} \label{Prop:overlap}
Let $T>0$ and $T_1 \in (0,T)$.
Let $v_j \in V_T$, $w_j \in L^2((0,T),\mathbb{R})$, $\Psi_j$ and $\xi_j$ be the associated solutions of  
(\ref{Schro_lin_eq}) and (\ref{Schro_O2_eq}) for $j=1,2$. We assume that $v_1$ is supported on $(0,T_1)$ and
$v_2$ is supported on $(T_1,T)$. Let $v:=v_1+v_2$, $w:=w_1+w_2$, $\Psi$ and $\xi$ be the associated solutions of  
(\ref{Schro_lin_eq}) and (\ref{Schro_O2_eq}). Then 
$\Psi(T)=0$ and $\xi(T)=\xi_1(T)+\xi_2(T)$.
\end{Prop}

\noindent \textbf{Proof of Proposition \ref{Prop:overlap}:}
We have $\Psi(T)=0$ because $v \in V_T$ (see (\ref{Psi_explicit}) and (\ref{def:VT})).
The control $v_1$ is supported on $(0,T_1)$ and belongs to $V_{T_1}$ thus
$\Psi_1$ is supported on $(0,T_1)\times(0,1)$ (see (\ref{Psi_explicit}) and (\ref{def:VT})).
The function $v_2$ is supported on $(T_1,T)$ thus $\Psi_2$ is supported on $(T_1,T)\times(0,1)$.
Therefore
$$(v_1+v_2)\mu(\Psi_1+\Psi_2)=v_1\mu\Psi_1+v_2\mu\Psi_2 \text{ on } (0,T)\times(0,1),$$
i.e. $\xi_1+\xi_2$ and $\xi$ solve the same Cauchy problem, thus $\xi=\xi_1+\xi_2$. \hfill  $\blacksquare$

\subsubsection{Proof of Proposition \ref{Prop:dir_perdue_co} in a simplified case}
\label{subsubsection_Ndir_simplified}

The strategy for the proof of Proposition \ref{Prop:dir_perdue_co}
is the same as in \cite{Crepeau-Cerpa}. It relies strongly on the rotation
of the lost directions, emphasized in Lemma \ref{lemme_sans_controle}. 
However, the strategy of (\cite{Crepeau-Cerpa}) needs to be adapted because there is no rotation phenomenon on our first lost direction.
In order to simplify the notations, we prove Proposition \ref{Prop:dir_perdue_co} in the case
$$N=2, \quad K_1=1, \quad K_2=2, \quad T_\sharp=2 T_{min}^2 + \frac{3 \pi}{\lambda_2 - \lambda_1},$$
where $T_{min}^2$ is defined in Lemma \ref{Lem:Q1<>0}. We will explain in Section \ref{subsubsection_Ndir} how it can be adapted for $N \geqslant 3$ and $K_1,...,K_N$ arbitrary. 
\\

\noindent 
Let $T$, $T_1$, $T_\theta$, $T_c$, $T_c^1>0$ be such that
\begin{equation} \label{hyp:T}
T > T_\sharp := 2 T_{min}^2 + \frac{3 \pi}{\lambda_2 - \lambda_1},
\end{equation}
\begin{equation} \label{hyp:T1}
\frac{\pi}{\lambda_2-\lambda_1} < T_1 < T  - \frac{2\pi}{\lambda_2 - \lambda_1} - 2 T_{min}^2,
\end{equation}
\begin{equation} \label{hyp:Tcth}
T_c < T_\theta, \quad T_c+T_\theta < \min\left\{ \frac{\pi}{\lambda_2-\lambda_1} ; T_1-\frac{\pi}{\lambda_2-\lambda_1} \right\}, 
\end{equation}
\begin{equation} \label{hyp:Tc1}
T_{min}^2 < T_c^1 < \frac{1}{2} \left( T-T_1-\frac{2\pi}{\lambda_2-\lambda_1} \right).
\end{equation}
Recall that $T_{min}^2$ is defined in Lemma \ref{Lem:Q1<>0}. Since $T_c^1 > T_{min}^2$, there exists controls $(\vpm,\wpm) \in L^2((0,T_c^1),\mathbb{R})^2$,
such that the associated solutions $\Psi^{\pm}$ and $\xi^{\pm}$ of (\ref{Schro_lin_eq}) and (\ref{Schro_O2_eq}) satisfy
\begin{equation} \label{ordre2_3_+-}
\begin{aligned}
&\Psi^{\pm}(T_c^1) = 0,
\\
&\lag \xi^{\pm}(T_c^1) , \psi_1(T_c^1) \rag = \pm i,
\\
&\lag \xi^{\pm}(T_c^1) , \psi_k(T_c^1) \rag =0, \quad  \forall k \geq 3.
\end{aligned}
\end{equation}
Indeed Lemma \ref{Lem:Q1<>0} implies the existence of $\vpm \in V_{T_c^1}$ such that $\widetilde{Q}_{1,T_c^1}^2(\vpm) = \pm 1$. Then, (\ref{espace_tangent_ordre2}) implies $\lag \xi^{\pm}(T_c^1) , \psi_1(T_c^1) \rag = \pm i$. Defining $\wpm$ as the solution of an adequate moment problem as in the proof of Proposition \ref{prop_dir_perdues} proves (\ref{ordre2_3_+-}).

\noindent
According to Proposition \ref{prop_dir_perdues}, 
there exists controls $(v^2,w^2) \in L^2((0,T_c),\mathbb{R})^2$ such that the associated solutions $(\Psi^2,\xi^2)$ 
of (\ref{Schro_lin_eq}) and (\ref{Schro_O2_eq}) satisfy
\begin{equation} \label{ordre2_3bis}
\begin{aligned}
&\Psi^2(T_c) = 0,
\\
&\lag \xi^2(T_c) , \psi_2(T_c) \rag \neq 0,
\\
&\lag \xi^2(T_c) , \psi_k(T_c) \rag =0, \quad  \forall k \geq 3.
\end{aligned}
\end{equation}

\noindent \emph{First step: Construction of a basis for $M^2=\text{Span}_{\mathbb{C}}(\psi_2(T))$, with nonoverlapping controls.}  
Let  
$$\begin{array}{ll}
\theta_1:=T-T_1, & \quad \theta_2:=T-T_1+\Tt, \\ 
\theta_3:=T-T_1+\frac{\pi}{\lambda_2-\lambda_1}, & \quad  \theta_4:=T-T_1+\Tt+\frac{\pi}{\lambda_2-\lambda_1}
\end{array}$$ 
and $(v^2_j,w^2_j):=(v^2_{\theta_j},w^2_{\theta_j})$ for $j=1,...,4$
with the notations of Lemma \ref{lemme_sans_controle} (in which $(T,\tilde{T})$ is replaced by $(T_c,T)$).
Then $\text{supp}(v^2_j) \subset (\theta_j,\theta_j+T_c)$ for $j=1,...,4$ and
$$T-T_1 = \theta_1 < \theta_1+T_c < \theta_2 < \theta_2+T_c < \theta_3 < \theta_3 +T_c <\theta_4 < \theta_4+T_c < T$$
(see (\ref{hyp:Tcth})), thus the supports do not overlapp:
\begin{equation} \label{overlap22}
\forall j_1 , j_2 \in \{1,2,3,4\} \text{ with } j_1 \neq j_2 \text{ then } 
\text{Supp}(v^2_{j_1}) \cap \text{Supp}(v^2_{j_2}) = \emptyset.
\end{equation}
We denote by $(\Psi^2_j,\xi^2_j)$ the associated solutions of (\ref{Schro_lin_eq}) and (\ref{Schro_O2_eq}).
Then, $\Psi^2_j(T)=0$ and $\xi_j^2(T)=\tilde{f}_j^2 + f_j^2$ for $j=1,...,4$ where (see Lemma \ref{lemme_sans_controle})
\begin{displaymath}
\begin{array}{ll}
\tilde{f}^2_1 = \lag \xi^2(T_c), \psi_1(T_c) \rag \psi_1(T), 
& \quad  
f^2_1 = e^{i(\lambda_2 - \lambda_1) (T-T_1)}\lag \xi^2(T_c) , \psi_2(T_c) \rag \psi_2(T) \neq 0,
\\
\tilde{f}^2_2 =  \tilde{f}^2_1,
& \quad
f^2_2 = e^{i (\lambda_2 - \lambda_1) \Tt} f^2_1,
\\
\tilde{f}^2_3 =  \tilde{f}^2_1,
& \quad
f^2_3 = e^{i (\lambda_2 - \lambda_1) \frac{\pi}{\lambda_2-\lambda_1} } f^2_1 = -f^2_1,
\\
\tilde{f}^2_4 =  \tilde{f}^2_1,
& \quad 
f^2_4 = e^{i (\lambda_2 - \lambda_1) (\frac{\pi}{\lambda_2-\lambda_1} + \Tt) } f^2_1 = -f^2_2.
\end{array}
\end{displaymath}
Moreover, (\ref{espace_tangent_ordre2}) imply that
\begin{equation} \label{Refjt}
\Re \lag \tilde{f}^2_j , \psi_1(T) \rag =
\Re \lag \xi^2(T_c) , \psi_1(T_c) \rag =
-\|\Psi^2(T_c)\|^2 =0, \quad \forall j=1,\dots ,4.
\end{equation}
Note that $(\lambda_2 - \lambda_1) \Tt \in (0, \pi)$, thus 
$(f^2_1,f^2_2)$ is a $\R$-basis of $M^2$. This leads to $M^2=  \bigcup_{j=1}^4  M^2_j$ where
\begin{equation} \label{def:M2j}
\begin{array}{l}
M^2_1 = \{ d^2_1 f^2_1 + d_2^2 f^2_2 ; d^2_1 \geq 0, d^2_2 \geq 0 \},
\\
M^2_2 = \{ d^2_1 f^2_2 + d_2^2 f^2_3 ; d^2_1>0, d^2_2 \geq 0 \},
\\
M^2_3 = \{ d^2_1 f^2_3 + d_2^2 f^2_4 ; d^2_1 \geq 0, d^2_2 \geq 0 \},
\\
M^2_4 = \{ d^2_1 f^2_4 + d_2^2 f^2_1 ; d^2_1>0, d^2_2 \geq 0 \}.
\end{array}
\end{equation}

\noindent \emph{Second step : Construction of a basis for $M^1$, with non overlapping controls.} 
The time interval $(T_c^1, T-T_1-T_c^1)$ has length
$(T-T_1-2T_c^1) >2\pi/(\lambda_2-\lambda_1)$ (see (\ref{hyp:Tc1})),
thus there exists an odd integer $k$ such that
\begin{equation} \label{def:AA}
\mathcal{T}:=\frac{k\pi}{\lambda_2-\lambda_1} \in (T_c^1, T-T_1-T_c^1).
\end{equation}
Let us consider the following controls
\begin{equation*}
(V_{\pm} , W_{\pm})(t) := 
\left\{ \begin{array}{ll}
(\vpm, \wpm)(t) & \quad \text{ if } t \in (0,T_c^1), \\
(0 ,0)  & \quad  \text{ if } t \in (T_c^1,\mathcal{T}), \\
(\vpm,\wpm)(t-\mathcal{T})  & \quad \text{ if } t \in (\mathcal{T}, \mathcal{T}+T_c^1),\\
(0 ,0)  & \quad \text{ if } t \in (\mathcal{T}+T_c^1,T),
\end{array} \right.
\end{equation*}
We denote by $(\Psi^1_{\pm},\xi^1_{\pm})$ the associated solutions of (\ref{Schro_lin_eq}) and (\ref{Schro_O2_eq}).
Then $\text{supp}(V_\pm) \subset [0,T-T_1)$ (see (\ref{def:AA})), thus
\begin{equation} \label{overlap12}
\forall j \in \{1,...,4\}, \text{Supp}(V_\pm) \cap \text{Supp}(v^2_j) = \emptyset.
\end{equation}
Then, $\Psi^1_\pm(T)=0$ and 
$$\xi^1_\pm(T)  = \pm 2 i \psi_1(T) + \langle \xi^\pm(T_c^1),\psi_2(T_c^1)\rangle [1+e^{i \mathcal{T}(\lambda_2-\lambda_1)}] \psi_2(T) 
= \pm 2 i \psi_1(T)$$
by Proposition \ref{Prop:overlap}, Lemma \ref{lemme_sans_controle} and (\ref{ordre2_3_+-}).

\noindent
As $M^1=i \text{Span}_{\mathbb{R}}(\psi_1(T))$, we can thus reach a $\mathbb{R}$-basis of $M^1$ with non-negative coefficients.
\\

\noindent \textit{Third step : Conclusion.} Let $z \in M$. We construct controls $(v,w) \in L^2((0,T),\mathbb{R})^2$
such that the associated solutions $(\Psi,\xi)$ of (\ref{Schro_lin_eq}) and (\ref{Schro_O2_eq}) satisfy
$\Psi(T)=0$ and $\xi(T)=z$. The proof relies  on the two following facts:
\begin{enumerate}
\item $\pm 2i\psi_1(T)$ and $f^2_j+\tilde{f}^2_j$ for $j=1,2,3,4$ are reachable states,
with controls such that their supports do not overlap (see (\ref{overlap22}) and (\ref{overlap12})),
\item any vector in $M$ is a linear combination of three of theses vectors, with only non negative coefficients
before $f^2_j+\tilde{f}^2_j$.
\end{enumerate}

There exists a unique $j \in \{1,2,3,4\}$ such that $z \in M^1 + M^2_j$ (see (\ref{def:M2j})).
Then, 
$$z=ix\psi_1(T)+d_1 f_j^2 + d_2 f_{j+1}^2 \text{ for some } d_1, d_2 \geqslant 0, x \in \mathbb{R}$$
with the convention $f^2_{5}=f^2_1$. We have
$$z=\Big( ix - d_1 \tilde{f}^2_j - d_2 \tilde{f}^2_{j+1} \Big)
+ d_1 \Big( \tilde{f}^2_j +  f_j^2 \Big) 
+ d_2 \Big( \tilde{f}^2_{j+1} +  f_{j+1}^2 \Big).$$  
As $\text{Re} ( \lag \tilde{f}^2_j , \psi_1(T) \rag) = 0$, for all $j=1,\dots,4$
(see (\ref{Refjt})), there exists $\kappa \in \{+,-\}$ and $c \geq 0$ such that
\begin{equation*}
 ix - d_1 \tilde{f}^2_j - d_2 \tilde{f}^2_{j+1} = \kappa 2 i  c \psi_1(T).
\end{equation*}
Then,
\begin{equation*}
z =  \kappa 2 i  c \psi_1(T) + d_1 ( f^2_j + \tilde{f}^2_j ) + d_2 ( f^2_{j+1} + \tilde{f}^2_{j+1} ),
\end{equation*}
i.e. $z$ is a linear combination of three states that are reachable with non overlapping controls. 
Hence the map
$$\Lambda_T(z):=(v,w):=
\Big( \sqrt{c} V_\kappa + \sqrt{d_1} v^2_j + \sqrt{d_2} v^2_{j+1} , 
c W_\kappa + d_1 w^2_j + d_2 w^2_{j+1} \big),$$
gives the conclusion. \hfill  $\blacksquare$
\\
\\

\subsubsection{Proof of Proposition \ref{Prop:dir_perdue_co} in the general case}
\label{subsubsection_Ndir}

Let us explain the adaptation of the strategy developed in Section \ref{subsubsection_Ndir_simplified} for $N\geq 3$. 
As previously, we denote by $K_1 < \dots < K_N$ the directions missed at the first order
and we explain how to reach a basis of missed directions on the second order (\ref{Schro_O2_eq}),
iteratively. In this proof, the term 'projection on $M^j$' denotes $\Im \lag \xi(T), \psi_1(T) \rag$ if $j=1$ and $K_1=1$ and $\lag \xi(T),\psi_{K_j}(T) \rag$ otherwise.
\newline

$\bullet$ The first step consists in reaching a $\R^+$ basis of $M^N$, 
the projections on $M^1, \dots, M^{N-1}$ being possibly non zero. 
This is done as in the first step of the proof of Proposition \ref{Prop:dir_perdue_co}, 
by designing four controls with non overlapping supports. 
It is done in any time $T_1 > \frac{\pi}{\lambda_{K_N}-\lambda_1}$. 
\\

$\bullet$ The $(k+1)^{th}$ step consists in reaching a $\R^+$ basis of $M^{N-k}$ 
while driving to zero the projections on $M^j$, for $j=N-k+1, \dots, N$. 
This can be done iteratively in the following way.
Let $(v^{(0)},w^{(0)})$ be as in Proposition \ref{prop_dir_perdues} for a sufficiently small time and for $j=N-k$. 
Then, the controls
\begin{equation*}
\big( v^{(1)} , w^{(1)} \big) := \big( v^{(0)} , w^{(0)} \big) + \big( v^{(0)}_{\theta} , w^{(0)}_{\theta} \big),
\quad \text{with } \theta = \frac{\pi}{\lambda_{K_N}-\lambda_1}
\end{equation*}
drive the projection on $M^N$ to zero while the projection on $M^{N-k}$ is still non zero.
This is ensured by Lemma \ref{lemme_sans_controle}. This is the same strategy as the second step of Section \ref{subsubsection_Ndir_simplified} where we drove the projection on $M^2$ to zero while the projection on $M^1$ was non zero. We iterate this construction
\begin{equation*}
\big( v^{(j+1)} , w^{(j+1)} \big) = \big( v^{(j)} , w^{(j)} \big) + \big( v^{(j)}_{\theta} , w^{(j)}_{\theta} \big),
\quad \text{with } \theta = \frac{\pi}{\lambda_{K_{N-j}}-\lambda_1}
\quad \text{ for } j=0,...,k-1.
\end{equation*}
Then the controls $\big( v,w \big) = \big( v^{(k)} , w^{(k)} \big)$ 
drive the projection on $M^N,\dots, M^{N-k+1}$ to zero 
while the projection on $M^{N-k}$ is still non zero. 
Finally, we can find $T_{\theta}$ sufficiently small such that $\big( v,w \big)$ and 
$\big( v_{T_{\theta}},w_{T_{\theta}} \big)$ have non overlapping supports 
and the four pairs of control 
$\big( v,w \big)$, $\big( v_{T_{\theta}},w_{T_{\theta}} \big)$, 
$\big( v_p , w_p \big)$ and $\big( v_{p+T_{\theta}}  ,w_{p+T_{\theta}} \big)$ 
with $p = \frac{\pi}{\lambda_{K_{N-k}}- \lambda_1}$ allow to conclude the $(k+1)^{th}$ step. 
This can be done in any time $T > \frac{\pi}{\lambda_{K_{N-k}}- \lambda_1} + \dots + \frac{\pi}{\lambda_{K_N}- \lambda_1}$.
\\

$\bullet$ The final step depends on the value of $K_1$. If $K_1 \geq 2$, we end with the same strategy. This step can be done in any time $T > \frac{\pi}{\lambda_{K_1}- \lambda_1} + \dots + \frac{\pi}{\lambda_{K_N}- \lambda_1}$ and leads to the expression (\ref{def_Tmin}) of $T_\sharp$.

\noindent
If $K_1 = 1$, the elementary brick of control cannot be designed in arbitrary small time (it was already the case in the second step of Section \ref{subsubsection_Ndir_simplified} where $(V_{\pm},W_{\pm})$ were constructed). In this case, the controls $(v^{(0)}, w^{(0)})$ have a time support greater than $T^2_{min}$. The iterative process then gives that this step can be done in any time
$T > 2^{N-1} T^2_{min} + \sum_{k=2}^N \frac{2^{k-2} \pi}{\lambda_{K_k}- \lambda_1}$ and leads to the expression (\ref{def_Tmin}) of $T_\sharp$.

\begin{figure}[H]
\centering
\vspace{-0.5cm}
\input{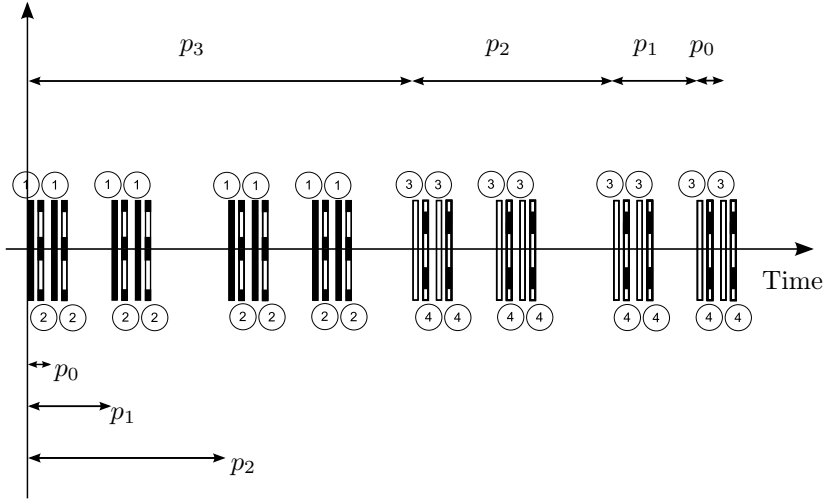}
\vspace{-1cm}
\caption{\label{fig:4dir} Support of controls when four directions are lost}
\end{figure}

%
%
Figure \ref{fig:4dir} illustrates the support of controls during the fourth step with $p_j := \frac{\pi}{\lambda_{K_{N-j}} - \lambda_1}$. The small rectangles indicates that the control is active. The phases of control associated with the same index define one of the four pair of controls $\big( v,w \big)$, $\big( v_{T_{\theta}},w_{T_{\theta}} \big)$, 
$\big( v_p , w_p \big)$ and $\big( v_{p+T_{\theta}}  ,w_{p+T_{\theta}} \big)$. The first rectangle (at the left) stands for the support of $(v^{(0)},w^{(0)})$.

\subsection{Proof of Theorem \ref{thm:control_Tlarge}}
\label{subsec:Fix_point}

Let $T>T_\sharp$ and $\psi_f \in H^3_{(0)}(0,1)$ be close enough to $\psi_1(T)$
(this will be precised later on). The goal of this section is the construction of 
$u \in L^2((0,T),\mathbb{R})$ such that
\begin{itemize}
\item the solution of (\ref{Schro:syst})-(\ref{IC}) satisfies $\psi(T)=\psi_f$,
\item $u$ tends to $0$ in $L^2((0,T),\mathbb{R})$ when $\psi_f \rightarrow \psi_1(T)$ in $H^3_{(0)}(0,1)$.
\end{itemize}
To simplify the notations, we assume $K_1=1$.
\\

Let $T_1 \in (T_\sharp,T)$ and $\delta_1>0$ associated to the map $\Gamma_{[T_1,T]}$ 
of Theorem \ref{thm:control_H}. From now on, we assume that
\begin{equation} \label{hyp:psif1}
\| \psi_f-\psi_1(T) \|_{H^3_{(0)}} < \delta_1.
\end{equation}
One may assume that $\delta_1$ is small enough so that condition (\ref{hyp:psif1}) implies
$\Re \langle \psi_f , \psi_1(T) \rangle >0$. We introduce the map
\begin{equation}
\label{def_Fpsi}
\begin{array}{|crcl}
F_{\psi_f}: & M \cap B_{L^2(0,1)}(0,\rho) & \rightarrow & M \\
            &                   z         & \mapsto     & \mathcal{P}_M[\psi_z(T)]
\end{array}
\end{equation}
where
\begin{itemize}
\item $\rho \in (0,1)$ will be chosen later on,
\item $\mathcal{P}_M:L^2(0,1) \rightarrow M$ is the $L^2$-orthogonal projection on $M$
$$\mathcal{P}_M(\zeta):=i \Im(\langle \zeta,\psi_1(T)\rangle)\psi_1(T) + 
\sum\limits_{j=2}^N \langle \zeta , \psi_{K_j}(T) \rangle \psi_{K_j}(T).$$
\item $\psi_z$ is the solution of  (\ref{Schro:syst})-(\ref{IC}) associated to the control $u_z$
defined by
$$u_z:=\left\lbrace \begin{array}{l}
\sqrt{\|z\|} v_z + \|z\| w_z \quad \text{ on } (0,T_1),\\
\Gamma_{[T_1,T]}( \psi_z(T_1) , \mathcal{P}_T[\psi_f] ) \quad \text{ on } (T_1,T).
\end{array}\right.$$
where $\|.\|$ is the $L^2(0,1)$-norm, $\Gamma_{[T_1,T]}$ is defined in Theorem \ref{thm:control_H}, $\mathcal{P}_T$ is defined by (\ref{def:P1T}) and
$$(v_z,w_z):=\Lambda_{T_1}\left( \frac{e^{iA(T-T_1)}z}{\|z\|} \right),$$
with $\Lambda_{T_1}$ defined in Proposition \ref{Prop:dir_perdue_co}.
\end{itemize}
Note that, for every $z$, $\mathcal{P}_T[\psi_z(T)]=\mathcal{P}_T[\psi_f]$.
Thus, our goal is to find $z_*$ such that $F_{\psi_f}(z_*)=\mathcal{P}_M[\psi_f]$.
\\

First, we check that the map $F_{\psi_f}$ is well defined when $\rho$ is small enough.

\begin{Prop} \label{Prop:FWD}
There exists $\rho>0$ such that, 
for every $\psi_f \in H^3_{(0)}(0,1)$ with (\ref{hyp:psif1}),
the map $F_{\psi_f}$ defined by (\ref{def_Fpsi}) is well defined and continuous on $M \cap B_{L^2(0,1)}(0,\rho)$.
\end{Prop}

\noindent \textbf{Proof of Proposition \ref{Prop:FWD}:} 
In order to prove that  $F_{\psi_f}$ is well defined, it is sufficient to find $\rho>0$ such that
\begin{equation} \label{FWD}
\|z\|<\rho \quad \Rightarrow \quad \|\psi_z(T_1)-\psi_1(T_1)\|_{H^3_{(0)}} < \delta_1.
\end{equation}
By Proposition \ref{WP-CYpb}, there exists $C_1, C_1'>0$ such that, for every $z \in M$,
$$\|\psi_z(T_1)-\psi_1(T_1)\|_{H^3_{(0)}} 
\leqslant C_1 \|u_z\|_{L^2(0,T_1)} 
\leqslant C_1' \sqrt{\|z\|}.$$
Thus, (\ref{FWD}) holds with $\rho:=\min\{1;(\delta_1/C_1')^2\}$. The continuity of $F_{\psi_f}$ is a consequence
of the continuity of $\Gamma_{[T_1,T]}$ and the continuity of the solutions of (\ref{Schro:syst})(\ref{IC_qlq})
with respect to the control $u$ and the initial condition $\psi_0$ (see (\ref{majo})). \hfill  $\blacksquare$
\\

One may assume $\rho$ small enough so that
$$\|z\|<\rho \quad  \Rightarrow \quad \Re \langle \psi_z(T),\psi_1(T) \rangle >0.$$
The goal of this section is the proof of the following result,
which proves Theorem \ref{thm:control_Tlarge}.

\begin{Prop} \label{Prop:z*}
There exists $\delta \in (0,\delta_1]$ such that, for every $\psi_f \in H^3_{(0)}(0,1)$ with
\begin{equation} \label{hyp:psif}
\| \psi_f-\psi_1(T) \|_{H^3_{(0)}} < \delta
\end{equation}
there exists $z_*=z_*(\psi_f) \in M \cap B_{L^2}(0,\rho)$ such that $F_{\psi_f}(z_*)=\mathcal{P}_M[\psi_f]$.
Moreover, $z_*(\psi_f) \rightarrow 0$ when $\psi_f \rightarrow \psi_1(T)$ in $H^3_{(0)}(0,1)$.
\end{Prop}

Combining $\mathcal{P}_T[\psi_z(T)]=\mathcal{P}_T[\psi_f]$, Proposition \ref{Prop:z*} and $||\psi_z(T)||_{L^2} = ||\psi_f||_{L^2}$ ends the proof of Theorem \ref{thm:control_Tlarge}.
The proof of  Proposition \ref{Prop:z*} requires the following preliminary result.

\begin{Prop} \label{Prop:z3/2}
There exists $\mathcal{C}>0$ such that,
for every $\psi_f \in H^3_{(0)}(0,1)$ with (\ref{hyp:psif1}) and $z \in M \cap B_{L^2(0,1)}(0,\rho)$,
we have
$$\| F_{\psi_f}(z) - z \| \leqslant \mathcal{C}[ \|\psi_f-\psi_1(T)\|_{H^3_{(0)}}^2 + \|z\|^{3/2} ].$$
\end{Prop}

\noindent \textbf{Proof of Proposition \ref{Prop:z3/2}:}\\

\noindent \emph{First step: Existence of $C_1>0$ such that 
\begin{equation} \label{z3/2:step1}
\| \psi_z(T_1) - \psi_1(T_1) - e^{iA(T-T_1)}z \|_{H^3_{(0)}} \leqslant C_1 \|z\|^{3/2}, 
\quad \forall z \in  M \cap B_{L^2(0,1)}(0,\rho).
\end{equation}}
Let $\Psi_z$, $\xi_z$ be the solution of (\ref{Schro_lin_eq}) and (\ref{Schro_O2_eq}) associated to the controls $v_z$ and $w_z$.
Then, $\Psi_z(T_1)=0$ and $\xi(T_1)=e^{iA(T-T_1)}z/\|z\|$. Explicit computations show that $\psi_z - \psi_1 - \sqrt{\|z\|} \Psi_z - \|z\| \xi_z$ is solution of (\ref{CY}) with control $u_z$, null initial condition and the following source term
$$(t,x) \mapsto ||z||^{3/2} w_z(t) \mu(x) \Psi_z(t,x) + ||z|| u_z(t) \mu(x) \xi_z(t,x).$$
In Proposition \ref{Prop:FWD}, $\rho$ was assumed to be smaller than $1$, thus Proposition \ref{WP-CYpb} implies that there exists $C>0$ such that
$$\| \psi_z - \psi_1 - \sqrt{\|z\|} \Psi_z - \|z\| \xi_z \|_{L^\infty((0,T),H^3_{(0)})} \leqslant C \|z\|^{3/2},
\forall z \in  M \cap B_{L^2(0,1)}(0,\rho).$$
which gives (\ref{z3/2:step1}).
\\

\noindent \emph{Second step: Existence of $C_2>0$ such that
\begin{equation} \label{z3/2:step2}
\|u_z\|_{L^2(T_1,T)} \leqslant C_2[ \|\psi_f-\psi_1(T)\|_{H^3_{(0)}} + \|z\| ],
\quad \forall z \in  M \cap B_{L^2(0,1)}(0,\rho).
\end{equation}}
The map $\Gamma_{[T_1,T]}$ is $C^1$ and $\Gamma_{[T_1,T]}(\psi_1(T_1),0)=0$, 
thus there exists $C>0$ such that, for every $ z \in  M \cap B_{L^2(0,1)}(0,\rho)$,
$$ \|u_z\|_{L^2(T_1,T)} = \| \Gamma_{[T_1,T]}( \psi_z(T_1) , \mathcal{P}_T[\psi_f] ) \|_{L^2(T_1,T)}
\leqslant C [ \|\psi_z(T_1)-\psi_1(T_1)\|_{H^3_{(0)}} + \| \mathcal{P}_T[\psi_f]\|_{H^3_{(0)}} ].$$
Explicit computations show that $\psi_z - \psi_1 - \sqrt{\|z\|} \Psi_z$ is solution of (\ref{CY}) with control $u_z$, null initial condition and the following source term
$$(t,x) \mapsto ||z|| w_z(t) \mu(x) \psi_1(t,x) + \sqrt{||z||} u_z(t) \mu(x) \Psi_z(t,x).$$
Thus Proposition \ref{WP-CYpb} implies that there exists $C>0$ such that
\begin{equation} \label{z3/2:T1}
\| \psi_z(T_1) - \psi_1(T_1) \|_{H^3_{(0)}} \leqslant \| \psi_z - \psi_1 - \sqrt{\|z\|} \Psi_z \|_{L^\infty((0,T),H^3_{(0)})} \leqslant C \|z\|.
\end{equation}
Then, $\PP_T[\psi_f]=\PP_T[\psi_f-\psi_1(T)]$ implies (\ref{z3/2:step2}).
\\

\noindent \emph{Third step: Existence of $C_3>0$ such that
\begin{equation} \label{z3/2:step3}
\| \psi_z - \psi_1 \|_{L^\infty((T_1,T),H^3_{(0)})} \leqslant
C_3 [ \|\psi_f-\psi_1(T)\|_{H^3_{(0)}} + \|z\| ],
\quad \forall z \in  M \cap B_{L^2(0,1)}(0,\rho).
\end{equation}}
Explicit computations show that $\psi_z - \psi_1$ is solution of (\ref{CY}) on $(T_1,T)$ with control $u_z$, initial condition $\psi_z(T_1) - \psi_1(T_1)$ and the following source term
$$(t,x) \mapsto u_z(t) \mu(x) \psi_1(t,x).$$
Using Proposition \ref{WP-CYpb} and (\ref{z3/2:T1}), we get a constant $C>0$ such that
$$\| \psi_z - \psi_1 \|_{L^\infty((T_1,T),H^3_{(0)})} \leqslant C [ \|z\| + \|u_z\|_{L^2(T_1,T)} ],
\forall z \in  M \cap B_{L^2(0,1)}(0,\rho),$$
which, together with (\ref{z3/2:step2}), gives (\ref{z3/2:step3}).
\\

\noindent \emph{Fourth step: Conclusion.}
Using the Duhamel formula, the commutativity between $e^{iAt}$ and $\mathcal{P}_M$ and
the isometry on $L^2(0,1)$ of $e^{iAt}$, we get for every $z \in  M \cap B_{L^2(0,1)}(0,\rho)$
\begin{align*}
\|F_{\psi_f}(z)-z\| & =  \| \mathcal{P}_M[\psi_z(T)] - z \|
\\ & \leqslant 
\| \mathcal{P}_M[e^{-iA(T-T_1)}\psi_z(T_1)] - z \| + \int_{T_1}^T |u_z(\tau)| \, \|\mathcal{P}_M[\mu \psi_z(\tau)] \| d\tau
\end{align*}
Then, using the relation $\mathcal{P}_M[\mu\psi_1(t)] \equiv 0$ 
(that holds because $\langle \mu \varphi_1,\varphi_{K_j}\rangle=0$ for $j=1,...,N$), Cauchy-Schwarz inequality and estimates (\ref{z3/2:step1}),(\ref{z3/2:step2}),(\ref{z3/2:step3}) it comes that
\begin{align*}
\|F_{\psi_f}(z)-z\| & \leqslant 
\| \mathcal{P}_M[\psi_z(T_1) - \psi_1(T_1)-e^{iA(T-T_1)}z ] \|
+ \int_{T_1}^T |u_z(\tau)| \, \|\mathcal{P}_M[\mu (\psi_z-\psi_1)(\tau)] \| d\tau
\\ & \leqslant 
\| \psi_z(T_1) - \psi_1(T_1) - e^{iA(T-T_1)}z \| 
\\ &  
+ \sqrt{T-T_1} \|u_z\|_{L^2(T_1,T)} \|\psi_z-\psi_1\|_{L^\infty((T_1,T),L^2)}
\\ & \leqslant 
C_1 \|z\|^{3/2} + \sqrt{T-T_1} C_2 C_3[ \|\psi_f-\psi_1(T)\|_{H^3_{(0)}} + \|z\| ]^2
\\ & \leqslant 
C(\rho) [ \|z\|^{3/2} + \|\psi_f-\psi_1(T)\|_{H^3_{(0)}}^2].
\end{align*}
This proves Proposition \ref{Prop:z3/2}. \hfill  $\blacksquare$
\\

\noindent \textbf{Proof of Proposition \ref{Prop:z*}:} We introduce the map
$$\begin{array}{|crcl}
G_{\psi_f}:  &  M \cap B_{L^2(0,1)}(0,\rho) & \rightarrow & M \\
              &           z                  & \mapsto     & z + \mathcal{P}_M[\psi_f]-F_{\psi_f}(z).
\end{array}$$
Our goal is to prove the existence of a fixed point $z_*=z_*(\psi_f)$ to the map $G_{\psi_f}$.
By Proposition \ref{Prop:z3/2}, there exists $\mathcal{C}>0$ (independent of $\psi_f$) such that, 
for every $z \in M \cap B_{L^2(0,1)}(0,\rho)$,
\begin{equation} \label{borne_G}
\begin{array}{ll}
\|G_{\psi_f}(z)\| 
& \leqslant \|z-F_{\psi_f}(z)\| + \|\mathcal{P}_M[\psi_f]\|
\\ & \leqslant 
\mathcal{C}[ \|\psi_f-\psi_1(T)\|_{H^3_{(0)}}^2 + \|z\|^{3/2} ]+ \|\psi_f-\psi_1(T)\|_{H^3_{(0)}}. 
\end{array}
\end{equation}
Let $\rho' \in (0,\rho)$ be such that 
\begin{equation} \label{hyp:rho'}
\mathcal{C} \sqrt{\rho'}<1/2
\end{equation}
and $\delta \in (0,\delta_1)$ be such that $\mathcal{C}\delta^2+\delta < \rho'/2$.
If $\psi_f$ satisfies (\ref{hyp:psif}), then $G_{\psi_f}$ maps continuously $M \cap B_{L^2(0,1)}(0,\rho')$ into itself.
The Brouwer fixed point theorem implies the existence of a fixed point
$z_*=z_*(\psi_f)$ of $G_{\psi_f}$ in $M \cap B_{L^2(0,1)}(0,\rho')$. We deduce from (\ref{borne_G}) and (\ref{hyp:rho'}) that
$$\|z_*(\psi_f)\| \leqslant 2[ \mathcal{C} \|\psi_f-\psi_1(T)\|_{H^3_{(0)}}^2 + \|\psi_f-\psi_1(T)\|_{H^3_{(0)}} ],$$
thus $z_*(\psi_f) \rightarrow 0$ when $\psi_f \rightarrow \psi_1(T)$ in $H^3_{(0)}(0,1)$. \hfill  $\blacksquare$

\section{Proof of Theorem \ref{thm:control_partiel_Tlarge_N}}
\label{sec:N}

In this section, we prove Theorem \ref{thm:control_partiel_Tlarge_N} when $\mu'(0)=\mu'(1)\neq 0$.
The case $\mu'(0)=-\mu'(1)\neq0$ may be proved similarly. The strategy is similar to the one of the previous
section, excepted that, for some lost directions, the second order may vanish and thus, we need to go
to a higher order. We prove that the third order is sufficient.

\subsection{Heuristic}

We consider a control $u$ of the form $u=\epsilon v+\epsilon^2 w+\epsilon^3 \nu$,
then, formally $\psi=\psi_1+\epsilon\Psi+\epsilon^2\xi+\epsilon^3 \zeta + o (\epsilon^3)$,
where $\Psi$ and $\xi$ solve (\ref{Schro_lin_eq}) and (\ref{Schro_O2_eq}) and
\begin{equation} \label{Schro_O3_eq}
\left\lbrace \begin{array}{ll}
i \partial_t \zeta = - \partial_x^2 \zeta - v(t) \mu(x) \xi - w(t) \mu(x) \Psi - \nu(t)\mu(x)\psi_1, &  (t,x) \in (0,T)\times(0,1), \\
\zeta(t,0)=\zeta(t,1)=0,                                                                             & t \in (0,T),                 \\
\zeta(0,x)=0.                                                                                        & x \in (0,1).          
\end{array}\right.
\end{equation}
We assume that $K \in \N^*$ satisfies $\lag \mu \varphi_1, \varphi_K \rag =0$ and that the quadratic form $Q_{K,T}^2$ vanishes on $V_T$ (see Proposition \ref{Prop:Q2<>0}).
Then, one may prove that for any $v \in V_T$, $\langle\zeta(T),\psi_K(T)\rangle=Q_{K,T}^3(v)$
where $Q_{K,T}^3$ is the cubic form (the index $3$ is related to the fact that $\zeta$ is the third order of the power series expansion)
$$Q_{K,T}^3(v):=\int_0^T v(t_1) \int_0^{t_1} v(t_2) \int_0^{t_2} v(t_3) h^3_{K,T}(t_1,t_2,t_3) dt_3 dt_2 dt_1,$$
$$h_{K,T}^3(t_1,t_2,t_3):=-i\sum\limits_{j_1=1}^\infty \sum\limits_{j_2=1}^\infty
B_{j_1,j_2} e^{i[(\lambda_K-\lambda_{j_1})t_1+(\lambda_{j_1}-\lambda_{j_2})t_2+(\lambda_{j_2}-\lambda_1)t_3]},$$
$$B_{j_1,j_2}:=\langle\mu\varphi_K,\varphi_{j_1}\rangle \langle\mu\varphi_{j_1},\varphi_{j_2}\rangle \langle\mu\varphi_{j_2},\varphi_1\rangle.$$

\begin{Prop} \label{Prop:Q2ouQ3}
Let $\mu \in H^3((0,1),\mathbb{R})$ and $K \in \mathbb{N}^*$ be such that $\mu'(0)=\mu'(1)\neq 0$
and $\langle \mu \varphi_1 , \varphi_K \rangle =0$.
Then,
\begin{itemize}
\item either, for every $N^*>0$, there exists $n\geqslant N^*$ such that
\begin{equation} \label{Q2ouQ3_1}
\langle \mu \varphi_{K},\varphi_n\rangle \langle\mu\varphi_n,\varphi_1\rangle \neq 0
\end{equation}
\item or, for every $N^*>0$, there exists $n_1,n_2\geqslant N^*$ such that
\begin{equation} \label{Q2ouQ3_2}
\langle \mu \varphi_{K},\varphi_{n_1}\rangle \langle\mu\varphi_{n_1},\varphi_{n_2}\rangle \langle\mu\varphi_{n_2},\varphi_1\rangle \neq 0.
\end{equation}
\end{itemize}
\end{Prop}

\noindent \textbf{Proof:} The proof relies on the equality (\ref{3IPP}).
If $K$ is odd, then (\ref{Q2ouQ3_1}) holds with $n$ odd and large enough. If $K$ is even and (\ref{Q2ouQ3_1}) does not hold, then (\ref{Q2ouQ3_2}) holds with $n_1$ odd, $n_2$ even, both large enough. \hfill  $\blacksquare$
\\

The previous and next propositions show that any lost direction (at the first order)
is recovered either at the second order, or at the third order.

\begin{Prop} \label{Prop:Q3<>0}
Let $\mu \in H^3((0,1),\mathbb{R})$, $K \in \mathbb{N}^*$ be such that
$$\langle \mu \varphi_{K},\varphi_{n_1}\rangle \langle\mu\varphi_{n_1},\varphi_{n_2}\rangle \langle\mu\varphi_{n_2},\varphi_1\rangle \neq 0
\text{ for some  } n_1, n_2 >K.$$ 
Then, $Q_{K,T}^3 \neq 0$ on $V_T$, $\forall T>0$.
\end{Prop}

\noindent \textbf{Proof of Proposition \ref{Prop:Q3<>0}:} To simplify the notations, we write $Q_T$ and $h$ instead of
$Q_{K,T}^3$ and $h^3_{K}$. 
Working by contradiction, we assume that $Q_{T} \equiv 0$ on $V_T$, for every $T<T^*$.
Then $\nabla Q_T (v) \perp V_T$, for every $v \in V_T$ and $T<T^*$.
Easy computations show that, for $v \in V_T$,
$$\nabla Q_{T}(v):t_3 \mapsto 
\int_{(0,T)^2} v(t_1) v(t_2) [ \tilde{h}(t_1,t_2,t_3)+\tilde{h}(t_1,t_3,t_2)+\tilde{h}(t_3,t_2,t_1)] dt_1 dt_2$$
where
$\tilde{h}(t_1,t_2,t_3):=h(t_1,t_2,t_3)1_{t_1>t_2>t_3}$. 
Let $v \in V_T$ with a compact support $(a,b) \subset (0,T)$.
For $t_3 \in (0,a)$, we have
\begin{equation} \label{grad_QT3}
\nabla Q_T(v)(t_3)= \sum\limits_{k_2=1}^\infty \alpha_{k_2}(v) e^{i(\lambda_{k_2}-\lambda_1)t_3}
\end{equation}
where
$$\alpha_{k_2}(v):=-i\int_{(a,b)^2} v(t_1) v(t_2)
\sum\limits_{k_1=1}^\infty 
B_{k_1,k_2} e^{i[(\lambda_K-\lambda_{k_1})t_1+(\lambda_{k_1}-\lambda_{k_2})t_2]} dt_2 dt_1.$$
We know that $\nabla Q_T(v)$ belongs to $\text{Adh}_{L^2(0,T)}(\text{Span}\{ e^{\pm i (\lambda_j-\lambda_1)t} ; j \in \JJ \})$
because $\nabla Q_T(v) \perp V_T$. The uniqueness of the decomposition on a Riesz basis ensures that (\ref{grad_QT3}) holds for all $t_3 \in (0,T)$.
For $t_3 \in (b,T)$, we have
$$\nabla Q_T(v)(t_3)=i\sum\limits_{k_1=1}^\infty \langle \mu \varphi_K,\varphi_{k_1}\rangle  Q^2_{k_1,T}(v)  e^{i(\lambda_K-\lambda_{k_1})t_3},$$
where $Q^2_{k_1,T}$ is defined in (\ref{def:Q2KT})-(\ref{def:h2K}).
Thus,
$$i\sum\limits_{k_1=1}^\infty  \langle \mu \varphi_K,\varphi_{k_1}\rangle Q^2_{k_1,T}(v) e^{i(\lambda_K-\lambda_{k_1})t_3}
=\sum\limits_{k_2=1}^\infty \alpha_{k_2}(v) e^{i(\lambda_{k_2}-\lambda_1)t_3},\; \forall b<t_3<T.$$
Notice that the frequencies $(\lambda_K-\lambda_{k_1})$ in the left hand side are negative when $k_1>K$,
and the frequencies $(\lambda_{k_2}-\lambda_1)$ in the right hand side are non-negative. Thus,
$$\langle \mu \varphi_K,\varphi_{k_1}\rangle Q^2_{k_1,T}(v)=0,  \forall k_1>K.$$
But $C^0_c(0,T) \cap V_T$ is dense in $V_T$, thus
\begin{equation} \label{absurd}
\langle \mu \varphi_K,\varphi_{k_1}\rangle Q^2_{k_1,T}\equiv 0 \text{ on } V_T,  \forall k_1>K.
\end{equation}

Let $n_1, n_2 > K$ be such that
$\langle \mu \varphi_{K},\varphi_{n_1}\rangle \langle\mu\varphi_{n_1},\varphi_{n_2}\rangle \langle\mu\varphi_{n_2},\varphi_1\rangle \neq 0$.
In particular, $\langle \mu \varphi_K,\varphi_{n_1}\rangle \neq 0$ 
and $Q^2_{n_1,T} \neq 0$ on $V_T$, for every $T>0$ by Proposition \ref{Prop:Q2<>0}.
This is in contradiction with (\ref{absurd}). \hfill  $\blacksquare$

\begin{rk}
Note that the third order may be necessary. For example, with
$\mu(x):=x-\langle x \varphi_1,\varphi_K\rangle \varphi_K/\varphi_1$,
where $K \in \mathbb{N}$ is even, we have
$\langle \mu \varphi_1,\varphi_n \rangle \langle \mu \varphi_n,\varphi_K \rangle =0$, $\forall n \in \mathbb{N}^*$,
thus $Q_{K,T}^2 \equiv 0$.
\end{rk}

\subsection{Reaching the missed directions at the second or third order}

We here only detail the changes with respect to the proof of Section \ref{sec:thm:control_Tlarge}. Using, (\ref{3IPP}) and the fact that $\mu'(0)=\mu'(1)\neq 0$, it comes that $\left\{ K \in \mathcal{N}_N \, ; \, \lag \mu \varphi_1 , \varphi_K \rag = 0 \right\}$ is finite. Thus, there exists $p \in \N^*$ and $K_1 < \dots < K_p \in \N^*$ such that for any $j \in \{1,\dots,p\}$, $\lag \mu \varphi_1, \varphi_{K_j} \rag = 0$. The estimate (\ref{3IPP}) also implies the existence of $C>0$ such that
\begin{equation*}
| \lag \mu \varphi_1 , \varphi_k \rag | \geq \frac{C}{k^3}, \quad \forall k \in \mathcal{N}_N- \{ K_1 , \dots , K_p \}.
\end{equation*}
For any $T>0$, Propositions \ref{Prop:Q2ouQ3} and \ref{Prop:Q3<>0} imply that for any $j \in \{1, \dots ,p\}$, if $Q_{K_j,T}^2$ vanishes on $V_T$, then $Q_{K_j,T}^3 \not \equiv 0$ on $V_T$.

\noindent
Let $\mathcal{K}_{(2)} := \left\{ j \in \{1, \dots,p\} \, ; \, Q_{K_j,T}^2 \not \equiv 0 \text{ on } V_T \right\}$ 
and $\mathcal{K}_{(3)} := \left\{1 , \dots, p \right\}-\mathcal{K}^{(2)}$.
The spaces $M^j$ and $M$ are defined as in (\ref{def:Mj}), (\ref{def:M}). Let us define
\begin{equation*}
M_{(2)} := \bigoplus_{j \in \mathcal{K}^{(2)} } M^j, 
\quad
M_{(3)} := \bigoplus_{j \in \mathcal{K}^{(3)} } M^j.
\end{equation*}
Thus, $M = M_{(2)} \oplus M_{(3)}$. 
Proposition \ref{Prop:dir_perdue_co} holds with $M$ replaced by $M_{(2)}$. 
The cubic form $Q_{K_j,T}^3$ satisfies $Q_{K_j,T}^3(-v) = -Q_{K_j,T}^3(v)$. Thus, one does not have to exploit the rotation phenomenon as in Proposition \ref{Prop:dir_perdue_co} and we can reach a basis with real non negative coefficients of $M^{(3)}$ on the third order in arbitrary time. More precisely, the following proposition holds.
\begin{Prop} \label{Prop:dir_perdue_O3_co}
Let $T>0$. There exists a continuous map 
$$\begin{array}{|cccl}
\tilde{\Lambda}_T: & M_{(3)} & \rightarrow & L^2((0,T),\mathbb{R})^2 \\
                   & z       & \mapsto     & (v,w,\nu)
\end{array}$$
such that, for every $z \in M_{(3)}$, the solutions $\Psi$, $\xi$ and $\zeta$ of
(\ref{Schro_lin_eq}), (\ref{Schro_O2_eq}) and (\ref{Schro_O3_eq}) satisfy $\Psi(T)=0$, $\xi(T)=0$ and $\zeta(T)=z$.
\end{Prop}
Finally, let us define the control $u_z$ by
$$u_z:=\left\lbrace \begin{array}{l}
\sqrt{\|z_2\|} v_{z_2} + \|z_2\| w_{z_2} + \|z_3\|^{1/3} \tilde{v}_{z_3} + \|z_3\|^{2/3} \tilde{w}_{z_3} + \|z_3\| \tilde{\nu}_{z_3}  \quad \text{ on } (0,T_1),\\
\Gamma_{[T_1,T]}( \psi_z(T_1) , \mathcal{P}_T[\psi_f] ) \quad \text{ on } (T_1,T).
\end{array}\right.$$
where $z_2 + z_3 = z$ with $(z_2,z_3) \in M_{(2)} \times M_{(3)}$ and
\begin{equation*}
(v_{z_2},w_{z_2}):=\Lambda_{T_1}\left( \frac{e^{iA(T-T_1)}z_2}{\|z_2\|} \right),\quad (\tilde{v}_{z_3},\tilde{w}_{z_3}, \tilde{\nu}_{z_3}):=\tilde{\Lambda}_{T_1}\left( \frac{e^{iA(T-T_1)}z_3}{\|z_3\|} \right).
\end{equation*}
Theorem \ref{thm:control_partiel_Tlarge_N} is then proved, as Theorem \ref{thm:control_Tlarge} in Section \ref{subsec:Fix_point}, using a fixed point argument. \hfill $\blacksquare$

\begin{rk}
As Proposition \ref{Prop:dir_perdue_co} is the only step requiring a minimal time, it has to be noticed that if $\mathcal{K}^{(2)} = \emptyset$, Theorem \ref{thm:control_partiel_Tlarge_N} holds in arbitrary time
\end{rk}

\section{A first step to the characterization of the minimal time}
\label{sec:caract}

In this section, we focus on the system
\begin{equation} \label{Schro:syst+S}
\left\lbrace \begin{array}{ll}
i \partial_t \psi(t,x) = - \partial_x^2 \psi(t,x) - u(t) \mu(x) \psi(t,x),  &  (t,x) \in \mathbb{R}\times(0,1),\\
\psi(t,0)=\psi(t,1)=0,                                                      & t \in \mathbb{R},                \\
s'(t)=u(t),                                                                 & t \in \mathbb{R},                \\
\end{array}\right.
\end{equation}
associated to the initial conditions
\begin{equation} \label{IC+S}
(\psi,s)(0)=(\varphi_1,0).
\end{equation}

We consider a dipolar moment $\mu \in H^3((0,1),\mathbb{R})$ such that
$\langle\mu\varphi_1,\varphi_1\rangle=0$ (for instance $\mu(x)=(x-1/2)$). We use the notation
$\mathcal{Q}_T$ instead of $\mathcal{Q}_{1,T}$ (see (\ref{def:QT(S)})-(\ref{def:kKT})),
$Q_T$ instead of $\widetilde{Q}^2_{1,T}$ (see (\ref{def:Q2KT})-(\ref{def:h2K})), 
$k(t,\tau)$ instead of $k_{1,T}(t,\tau)$
and the spaces
$$V_T^1 :=\left\{ v \in L^2(0,T) ; \int_0^T v(t) e^{i \omega_j t} dt = 0, \forall j \in \JJ \cup\{1\} \right\}$$
$$\mathcal{V}_T:=\left\{ S \in L^2((0,T),\mathbb{R}) ; 
\int_0^T S(t)e^{i\omega_j t} dt =0, \forall j \in \JJ  \right\}$$
where $\JJ$ is defined by (\ref{def:J}).
We introduce the quantities
$$\tilde{T}_{min}^1:=\sup\{ T \geqslant 0 ; \mathcal{Q}_{T} \leqslant 0 \text{ on } \mathcal{V}_T \},$$
\begin{equation} \label{def:Tmin2}
\tilde{T}_{min}^2:=\inf\{ T \geqslant 0 ; \exists S_{\pm} \in \mathcal{V}_T \cap H^1_0(0,T) \text{ such that } 
\mathcal{Q}_{T}(S_\pm)=\pm 1 \}.
\end{equation}
Lemma \ref{Lem:T_2} ensures that $\tilde{T}_{min}^1>0$ and the following proposition justifies the existence of $\tilde{T}_{min}^2$.

\begin{Prop} \label{Prop:Tmin12_comp}
Let $\mu \in H^3((0,1),\mathbb{R})$ be such that $\langle\mu\varphi_1,\varphi_1\rangle=0$.
For every $T>2/\pi$, there exists $S_\pm \in \mathcal{V}_T \cap H^1_0(0,T)$ such that $\mathcal{Q}_{T}(S_\pm)=\pm 1$;
or, equivalently, there exists $v_\pm \in V_T^1$ such that $Q_{T}(v_\pm)=\pm 1$.
Thus,
$$0<T_1^*<\tilde{T}_{min}^1 \leqslant \tilde{T}_{min}^2 \leqslant \frac{2}{\pi},$$
where $T_1^*$ was defined in Lemma \ref{Lem:T_2}.
\end{Prop}

This proposition may be proved as Lemma \ref{Lem:Q1<>0}.
The goal of this section is the proof of the following theorem.

\begin{thm} \label{thm:caract}
Let $\mu \in H^3((0,1),\mathbb{R})$ be such that 
\begin{equation} \label{hyp_mu_caract}
\langle\mu\varphi_1,\varphi_1\rangle=0
\quad \text{ and } \quad
\exists c>0 \text{ such that } \frac{c}{k^3} \leqslant |\langle\mu\varphi_1,\varphi_k\rangle|, \forall k \in \JJ.
\end{equation}
\begin{enumerate}
\item For every $T<\tilde{T}_{min}^1$, there exists $\epsilon>0$ such that, for every $u \in L^2((0,T),\mathbb{R})$ with (\ref{upt})
the solution of (\ref{Schro:syst+S})-(\ref{IC+S}) satisfies $(\psi,s)(T) \neq ( [\sqrt{1-\delta^2} + i \delta ] \psi_1(T) , 0)$
for every $\delta>0$.
\item If, moreover $\JJ=\mathbb{N}^*-\{1\}$, then, 
for every $T>\tilde{T}_{min}^2$, the system (\ref{Schro:syst+S}) is controllable in $H^3_{(0)}(0,1) \times \mathbb{R}$,
locally around the ground state $(\psi=\psi_1,s\equiv 0)$,
in time $T$, with controls $u \in L^2((0,T),\mathbb{R})$.
\end{enumerate}
In particular, when $\JJ=\mathbb{N}^*-\{1\}$, the minimal time $T_{min}$ required for the local controllability satisfies
$T_{min} \in [\tilde{T}_{min}^1,\tilde{T}_{min}^2]$.
\end{thm}

\begin{rk}
The equality between $\tilde{T}_{min}^1$ and $\tilde{T}_{min}^2$
is an open problem, equivalent to the question addressed in the next paragraph.

Let $P_T$ be the orthogonal projection from $L^2((0,T),\mathbb{R})$ to 
the closed subspace $\mathcal{V}_T$ and $\mathcal{K}_T$ be the compact self adjoint operator on $L^2((0,T),\mathbb{R})$ defined by
$$\mathcal{K}_T:=P_T \left[ t \mapsto \int_0^t k(t,\tau) S(\tau) d\tau \right].$$
Recall that $A_1$ is defined by (\ref{hyp_mu_0}). We know that 
\begin{itemize}
\item for any $T<\tilde{T}_{min}^1$ all the eigenvalues of $\mathcal{K}_T$ are $<A_1$
(see the first statement of Theorem \ref{thm:caract}),
\item for any $T>\tilde{T}_{min}^1$, the largest eigenvalue of $\mathcal{K}_T$ is $>A_1$.
(by definition of $\tilde{T}_{min}^1$).
\end{itemize}
For $T>\tilde{T}_{min}^1$, does the associated eigenvector belong to $H^1_0((0,T),\mathbb{R})$?
\end{rk}

The proof of the second statement of Theorem \ref{thm:caract} may be done exactly as the proof of Theorem \ref{thm:control_Tlarge} in Section \ref{sec:thm:control_Tlarge}.
Indeed, when $\JJ=\mathbb{N}^*-\{1\}$, then
\begin{enumerate}
\item the vector space $M$ of lost directions (at the first order) is $i \mathbb{R} \psi_1(T)$,
\item for any $T_1 \in (\tilde{T}^2_{min},T)$, the controls $S_{\pm} \in \mathcal{V}_{T_1} \cap H^1_0(0,T_1)$
allow to reach the states $\pm i \psi_1(T_1)$ with the second order term;
moreover, $(i\psi_1(T_1),-i\psi_1(T_1))$ is an '$\mathbb{R}^+$-basis' of $M$.
\end{enumerate}
Thus, in this section, we focus only on the proof of the first statement of Theorem \ref{thm:caract},
which is a direct consequence of the following result.

\begin{thm} \label{Thm:aux_bis}
Let $\mu \in H^3((0,1),\mathbb{R})$ that satisfies (\ref{hyp_mu_caract}).
For every $T<\tilde{T}_{min}^1$, 
there exists $\epsilon>0$ such that
for every $s \in H^1((0,T),\mathbb{R})$ with $s(0)=0$ and $\| \dot{s} \|_{L^2} < \epsilon$,
the solution of the Cauchy problem (\ref{Schro:syst_aux}) satisfies
$\widetilde{\psi}(T) \neq ( \sqrt{1-\delta^2} + i \delta ) \psi_1(T)$, $\forall \delta>0$.
\end{thm}

In section \ref{subsec:prel}, we state a preliminary result for the proof of Theorem \ref{Thm:aux_bis},
which is detailled in section \ref{subsec:thm10}.

\subsection{Preliminaries}
\label{subsec:prel}

For $T>0$ and $\eta>0$, we introduce the sets
\begin{equation} \label{def:VTeta}
\VV_{T,\eta}:=\left\{ S \in L^2(0,T); 
\left\| \left( \int_0^T S(t) e^{i \omega_j t} dt \right)_{j \in \JJ} \right\|_{l^2} \leqslant \eta \|S\|_{L^2(0,T)}
\right\}
\end{equation}
where $\JJ$ is defined in (\ref{def:J}). 

\begin{Prop} \label{Prop:Coercivite_T<Tmin1}
For every $T<\tilde{T}_{min}^1$, there exists $\lambda=\lambda(T), \eta=\eta(T)>0$ such that
\begin{equation} \label{coer_VT}
\mathcal{Q}_{T}(S) \leqslant -\lambda(T) \|S\|_{L^2(0,T)}^2, \quad \forall S \in \mathcal{V}_T,
\end{equation}
\begin{equation} \label{coer_VTeta}
\mathcal{Q}_{T}(S) \leqslant -\frac{\lambda(T)}{2} \|S\|_{L^2(0,T)}^2, \quad \forall S \in \mathcal{V}_{T,\eta}.
\end{equation}
\end{Prop}

This proposition may be proved with the formalism of Legendre quadratic forms (see \cite{Bonnans_book}).
For this article to be self contained, we propose an elementary proof in Appendix \ref{appendix_coercivite}.
\\

\subsection{Proof of Theorem \ref{Thm:aux_bis}}
\label{subsec:thm10}

Let $T<\tilde{T}_{min}^1$. We proceed as in the proof of Theorem \ref{Thm:aux}.
Working by contradiction, we assume that,  for every $\epsilon>0$, 
there exists $s_\epsilon \in H^1(0,T)$ with $s_\epsilon(0)=0$ and $\|\dot{s}_\epsilon\|_{L^2}<\epsilon$ 
such that the solution $\widetilde{\psi}_{\epsilon}$ of  (\ref{Schro:syst_aux}) satisfies 
\begin{equation} \label{target_imp_aux}
\widetilde{\psi}_\epsilon(T) = ( \sqrt{1-\delta_\epsilon^2} + i \delta_\epsilon ) \psi_1 (T), 
\end{equation}
for some $\delta_\epsilon>0$.
\\

\noindent \emph{First step : For $\epsilon>0$ small enough, $s_\epsilon \in \mathcal{V}_{T,\eta}$
(with $\eta=\eta(T)$ as in Proposition \ref{Prop:Coercivite_T<Tmin1}).}
Using (\ref{hyp_mu_caract}), Proposition \ref{Prop:approx_Comp1k} and (\ref{target_imp_aux}) we have
\begin{align*}
\left\| \left( \int_0^T s_\epsilon(t) e^{i \omega_j t} dt \right)_{j \in \JJ} \right\|_{l^2}
& \leqslant
C \left\| \left( \omega_j \langle \mu \varphi_1,\varphi_j \rangle \int_0^T s_\epsilon(t) e^{i \omega_j t} dt \right)_{j \in \JJ} \right\|_{h^1}
\\ & \leqslant
C \left\| \left(  \langle \widetilde{\psi}_\epsilon(T) , \psi_j(T) \rangle \right)_{j \in \JJ}  \right\|_{h^1}
+ \underset{\epsilon \to 0}{o}(\|s_\epsilon\|_{L^2})
\\ & = \underset{\epsilon \to 0}{o}(\|s_\epsilon\|_{L^2})
\label{s_dans_VT}
\tag{\theequation} \addtocounter{equation}{1}
\end{align*}
which gives the conclusion.
\\

\noindent
\emph{Second step : Conclusion.}
Using (\ref{approx_comp1_O2}) with $K=1$, the first step and (\ref{coer_VTeta}) it comes that
\begin{align*}
0 < \delta_{\epsilon} &= \Im \langle \widetilde{\psi}_{\epsilon}(T), \psi_1(T) \rangle
\\
&= \mathcal{Q}_T(s_{\epsilon}) + \underset{\epsilon\rightarrow 0}{o}( \|s_\epsilon\|_{L^2}^2 )
\\ 
& \leqslant - \frac{\lambda(T)}{2} \|s_\epsilon\|_{L^2}^2 + \underset{\epsilon\rightarrow 0}{o}( \|s_\epsilon\|_{L^2}^2),
\end{align*}
which is impossible for $\epsilon$ small enough. This ends the proof of Theorem \ref{Thm:aux_bis}.
 \hfill  $\blacksquare$
\\

\subsection{Comments about generalizations}

Let us consider a situation in which the first order misses exactly $N\geq2$ directions associated to the indexes $K_1, \dots, K_N$. Let $Q_{K_1,T}^2,\dots, Q_{K_N,T}^2$ be the associated complex-valued quadratic forms. A natural candidate for the minimal time $T_{min}$ could be the minimal time $\tilde{T}_{min}$ for the image of 
$$ \big( Q_{K_1,T}^2,\dots, Q_{K_N,T}^2 \big) : V_T \to \C^N$$
to cover $\C^N$. The positive controllability result in time $T>\tilde{T}_{min}$ could be proved with the technics of this article. The negative controllability result in time $T<\tilde{T}_{min}$ is more difficult. To transfer an impossible motion from the second order to the nonlinear system we need a coercivity property which is not obvious in this case.

\section{Conclusion, open problems, perspectives}
\label{sec:ccl}

In Theorem \ref{Main_thm_0}, we have proposed a general context for 
the local controllability of the system (\ref{Schro:syst}) 
to require a positive minimal time. 
This statement extends Coron's previous result in  \cite{JMC-CRAS-Tmin} because:
\begin{enumerate}
\item it does not use the variables $(s,d)$ in the state,
\item the control $u$ has to be small in $L^2$ (not in $L^\infty$),
\item $\mu(x)$ is not necessarily $(x-1/2)$.
\end{enumerate}
The validity of the conclusion without the assumption $A_K \neq 0$ is an open problem.
\\

In Theorem \ref{sec:thm:control_Tlarge}, we have proposed a sufficient condition
for the system (\ref{Schro:syst}) to be controllable around the ground state in large time.
This sufficient condition is compatible with the general context of  Theorem \ref{Main_thm_0},
thus there exists a large class of functions $\mu$ for which local controllability
holds in large time, but not in small time.
\\

The existence of a positive minimal time for the controllability is closely related to a second 
order approximation of the solution. When a direction is not controllable
neither at the first order, nor at the second one, then it is recovered at the third one,
and no minimal time is required.
\\

The characterization of the minimal time for the local controllability around the ground state
is essentially an open problem.
A first step has been done in this article, when only the first direction is lost.
\\

In \cite{Crepeau-Cerpa}, Crépeau and Cerpa prove the local controllability of the KdV equation,
with boundary control. When the length of the domain is critical, the linearized system is not controllable
along a finite number of directions, but all of them are recovered at the second or third order.
The existence of a positive minimal time, required for the local controllability is an open problem.
The technics developed in this article may be helpful for this question.

\begin{center}
\textbf{Acknowledgments}
\end{center}

The authors thank Jean-Michel Coron for having attracted their attention to this problem.

\appendix

\section{Trigonometric moment problems}

In this article, we use several times the following result 
(see, for instance \cite[Corollary 1 in Appendix B]{KB-CL} for a proof).

\begin{Prop} \label{Cor:haraux1}
Let $(\omega_{k})_{k \in \mathbb{N}^*}$ be an increasing sequence of $[0,+\infty)$ such that 
$\omega_{k+1} - \omega_{k} \rightarrow + \infty$ when $k \rightarrow + \infty$ and $\omega_1=0$.
Let $l^2_r(\mathbb{N}^*,\mathbb{C}):=\{ d=(d_k)_{k \in \mathbb{N}^*} \in l^2(\mathbb{N}^*,\mathbb{C}) ; d_1 \in \mathbb{R} \}$.
\begin{enumerate}
\item For every $T>0$, there exists a continuous linear map
$$\begin{array}{cccc}
L_T: & l^2_r(\mathbb{N}^*,\mathbb{C}) & \rightarrow & L^2((0,T),\mathbb{R}) \\
     &               d                & \mapsto     & L_T(d)
\end{array}$$
such that, for every $d=(d_k)_{k \in \mathbb{N}^*} \in l^2(\mathbb{N}^*,\mathbb{C})$,
the function $v:=L_T(d)$ solves
$$\int_0^T v(t) e^{i \omega_k t} dt = d_k, \forall k \in \mathbb{N}^*.$$
\item For every $T>0$ there exists a constant $C=C(T)$ such that (Ingham inequality)
$$\sum\limits_{k=1}^\infty |a_k|^2 \leqslant 
C \int_0^T \left| \sum\limits_{k=1}^\infty a_k e^{i\omega_k t} \right|^2 dt,
\quad \forall (a_k)_{k \in \mathbb{N}^*} \in l^2(\mathbb{N}^*,\mathbb{C}).$$
\\
\end{enumerate}
\end{Prop}

\section{Proof of Lemma \ref{lemme:KB-CL}}
\label{appendix_lemme_KB-CL}

This appendix is devoted to the proof of Lemma \ref{lemme:KB-CL}. It is a straightforward adaptation of \cite[Lemma 1]{KB-CL}. By definition, 
\begin{equation*}
F(t) = \sum_{k=1}^{\infty} \Big( \int_0^t \langle f(\tau), \varphi_k \rangle e^{i \lambda_k \tau} d \tau \Big) \varphi_k, \quad \text{in } L^2(0,1).
\end{equation*}
For almost every $\tau \in (0,T)$, $f(\tau) \in H^1$ and
\begin{align*}
\langle f(\tau) , \varphi_k \rangle &= \sqrt{2} \int_0^1 f(\tau,x) \sin(k \pi x) d x
\\
&= \frac{-\sqrt{2}}{k \pi} \big( (-1)^k f(\tau,1) - f(\tau,0) \big) + \frac{\sqrt{2}}{k \pi} \int_0^1 f'(\tau,x) \cos(k \pi x) d x.
\end{align*}
Thus,
\begin{align*}
||F(t)||_{H^1_0} &= \Big| \Big| \int_0^t \langle f(\tau) ,\varphi_k \rangle e^{i \lambda_k \tau} d \tau \Big| \Big|_{h^1}
\\
& \leqslant \frac{\sqrt{2}}{\pi} \Big( \Big| \Big| \int_0^t f(\tau,1) e^{i \lambda_k \tau} d \tau \Big| \Big|_{\ell^2} + \Big| \Big| \int_0^t f(\tau,0) e^{i \lambda_k \tau} d \tau \Big| \Big|_{\ell^2} \Big)
\\
&+ \frac{1}{\pi} \Big| \Big| \int_0^t \langle f'(\tau) , \sqrt{2} \cos(k \pi x) \rangle e^{i \lambda_k \tau} d \tau \Big| \Big|_{\ell^2}.
\end{align*}
As $(\sqrt{2} \cos(k \pi x))_{k\in \N^*}$ is orthonormal in $L^2(0,1)$,
\begin{align*}
\Big| \Big| \int_0^t \langle f'(\tau) , \sqrt{2} \cos(k \pi x) \rangle e^{i \lambda_k \tau} d \tau \Big| \Big|_{\ell^2} 
&= \left( \sum_{k=1}^{\infty} \Big| \int_0^t \langle f'(\tau) , \sqrt{2} \cos(k \pi x) \rangle e^{i \lambda_k \tau} d \tau \Big|^2 \right)^{1/2}
\\
&\leqslant \left( \sum_{k=1}^{\infty} t  \int_0^t \big| \langle f'(\tau) , \sqrt{2} \cos(k \pi x) \rangle \big|^2 d \tau \right)^{1/2}
\\
&\leqslant \sqrt{t} \left( \int_0^t ||f'(\tau)||_{L^2}^2 d \tau \right)^{1/2}
\\
&\leqslant \sqrt{t} || f ||_{L^2((0,t),H^1)}.
\end{align*}
Finally \cite[Appendix B, Corollary 4]{KB-CL} imply
\begin{align*}
||F(t)||_{H^1_0} &\leqslant \frac{\sqrt{2} C(t)}{\pi} \big( ||f'( \cdot, 1)||_{L^2(0,t)} + ||f'(\cdot, 0)||_{L^2(0,t)} \big) + \frac{\sqrt{2}}{\pi} ||f||_{L^2((0,t),H^1)}
\\
&\leqslant c_1(t) ||f||_{L^2((0,t),H^1)}
\end{align*}
where $c_1(t)$ is bounded for $t$ lying in bounded intervals. This proves that $F(t) \in H^1_0(0,1)$ for every $t \in [0,T]$ and that $t \mapsto F(t) \in H^1_0$ is continuous at $t=0$. The  continuity at any $t \in [0,T]$ may be proved similarly.
\hfill $\blacksquare$

\section{Proof of Proposition \ref{Prop:Coercivite_T<Tmin1}}
\label{appendix_coercivite}

This appendix is devoted to the proof of Proposition \ref{Prop:Coercivite_T<Tmin1}. The proof is divided in two steps. First, using a maximizing sequence we prove (\ref{coer_VT}). Then, solving an adequate moment problem, we prove (\ref{coer_VTeta}).
\\

\noindent \emph{First step: Proof of (\ref{coer_VT}).} For $T \in (0,\tilde{T}_{min}^1)$, we define the quantity $\lambda(T) \geqslant 0$ by
\begin{equation}
\label{def_lambda}
- \lambda(T) := \sup \{ \QQ_T(S) ; S \in \VV_T , ||S||_{L^2(0,T)}=1 \}.
\end{equation}

First, let us emphasize that,
if $ \lambda(T) \leqslant 0$, then, there exists $S \in \VV_{T}$ such that $||S||_{L^2(0,T)}=1$ and $\QQ_T(S) = \lambda(T)$
(consider a weak $L^2(0,T)$-limit, of a maximizing sequence and use the compactness of the operator 
$K:L^2(0,T)\rightarrow L^2(0,T)$ defined by
$KS : t \mapsto  \int_0^t S(\tau) k(t,\tau) d \tau$).

Let us assume that there exists $T \in (0,\tilde{T}_{min}^1)$ such that $\lambda(T)=0$.
Let $T_1 \in (T,\tilde{T}_{min}^1)$. 
Let $S_* \in \VV_{T}$ such that $\|S_*\|_{L^2(0,T)}=1$ and $\QQ_T(S_*) = 0$.
We extend $S_*$ on $(T,T_1)$ by zero.
Then, $S_* \in \VV_{T_1}$ and $\QQ_{T_1}(S_*)= \max \{\QQ_{T_1}(S);S\in\VV_T\}=0$ thus
(Euler equation) $\nabla \QQ_{T_1} (S_*) \perp \VV_{T_1}$, i.e. there exists a unique
sequence $(a_j)_{j \in \JJ-\{1\}} \in l^2$ such that
$$\nabla \QQ_{T_1} S_* (t) = \sum\limits_{j \in \JJ-\{1\}} a_j e^{i \omega_j t} \text{ in } L^2(0,T_1).$$
However, we have
$$\nabla \QQ_{T_1} (S_*) (t) = -A_1 S_*(t)+\int_0^t S_*(\tau)k(t,\tau)d\tau, \forall t \in (0,T_1).$$
In particular, $\nabla \QQ_{T_1} (S_*) \equiv 0$ on $(T,T_1)$ thus (Ingham inequality, see Proposition \ref{Cor:haraux1}) $a_j \equiv 0$.
We have proved that
$$S_*(t)=\frac{1}{A_1}\int_0^t S_*(\tau)k(t,\tau)d\tau, \forall t \in (0,T).$$
Thus, $S_*(0)=0$, $S_* \in H^1((0,T),\mathbb{R})$ and $S_*'$ satisfies the same relation.
Iterating this result, we get $S_*^{(n)}(0)=0$ and $S_*^{(n)} \in \text{Ker}(-A_1 Id + K)$ for every $n \in \mathbb{N}$.
But $K$ is compact, so $\text{dim}[\text{Ker}(-A_1 Id + K)] < +\infty$. Thus 
there exists $N \in \mathbb{N}^*$ and $a_0,\cdots,a_{N-1} \in \mathbb{R}$ such that 
\begin{equation*}
\left\{
\begin{aligned}
S_*^{(N)} &= a_0 S_* + a_1 S_*' + \cdots + a_{N-1} S_*^{(N-1)}
\\
S_*(0)& =0 , \cdots , S_*^{(N-1)}(0) = 0
\end{aligned}
\right.
\end{equation*}
\\
Therefore $S_*=0$, which is a contradiction.
\\

\noindent \emph{Second step: Proof of (\ref{coer_VTeta}):}
Let $\eta>0$ and $S \in \VV_{T,\eta}$ with $\|S\|_{L^2}=1$. 
Let $d:=(d_k)_{k \geqslant 2}$ be defined by
$$d_k:=\int_0^T S(t) e^{i\omega_k t} dt, \forall k \geqslant 2.$$
Then $\|d\|_{l^2} \leqslant \eta$. 
Let $\tilde{S}:=L_T(d)$ and $S_0:=S-\tilde{S}$, where $L_T$ is as in Proposition \ref{Cor:haraux1}.
Let $C(T):=\|L_T\|$. We have
\begin{equation} \label{bornes}
\| \tilde{S} \|_{L^2} \leqslant C(T) \eta
\quad \text{ and } \quad
1-C(T) \eta \leqslant \| S_0 \|_{L^2} \leqslant 1+C(T)\eta. 
\end{equation}
Using the first step and Cauchy-Schwarz inequality, we get
$$\begin{array}{ll}
\QQ_T(S) 
& = \QQ_T(S_0+\tilde{S})
\\ & =
\QQ_T(S_0)+\QQ_T(\tilde{S})+  \int_0^T S_0(t) \int_0^t \tilde{S}(s) k(t,s) ds dt + 
\int_0^T \tilde{S}(t) \int_0^t S_0(s) k(t,s) ds dt
\\ & \leqslant
- \lambda(T) \|S_0\|_{L^2}^2 + \frac{T}{2} \|k\|_\infty \|\tilde{S}\|_{L^2}^2 + 2 T \|k\|_{\infty} \|S_0\|_{L^2} \|\tilde{S}\|_{L^2}
\\ & \leqslant
- \lambda(T) [1-C(T) \eta]^2 + \frac{T}{2} \|k\|_\infty C(T)^2 \eta^2 + 2 T \|k\|_{\infty} [1+C(T)\eta] C(T)\eta.
\end{array}$$
Thus, for $\eta$ small enough, we get $\QQ_T(S) \leqslant -\lambda(T)/2 < 0$.
\hfill  $\blacksquare$

\bibliography{biblio3}
\bibliographystyle{plain}


\end{document}